\RequirePackage{fix-cm}
\documentclass[smallextended]{svjour3}       % onecolumn (second format)
\smartqed  % flush right qed marks, e.g. at end of proof
\usepackage{graphicx}
\usepackage{xcolor}
\usepackage{amsmath,amssymb,cite,url,booktabs}
\usepackage[numbers]{natbib}
\newcommand{\reals}{{\mbox{\bf R}}}
\newcommand{\BEQ}{\begin{equation}}
\newcommand{\EEQ}{\end{equation}}
\newcommand{\BIT}{\begin{itemize}}
\newcommand{\EIT}{\end{itemize}}
\newcommand{\ie}{{\it i.e.}}
\newcommand{\eg}{{\it e.g.}}
\newcommand{\Tr}{\mathop{\bf Tr}}
\newcommand{\argmin}{\mathop{\rm argmin}}
\newcommand{\Co}{{\mathop {\bf Co}}}
\newcommand{\ones}{\mathbf 1}
\newcommand{\Expect}{\mathop{\bf E{}}}
\newcommand{\Prob}{\mathop{\bf Prob}}
\newcommand{\VaR}{\mathop{\bf VaR}}
\newcommand{\CVaR}{\mathop{\bf CVaR}}
\newcommand{\EVaR}{\mathop{\bf EVaR}}

\newcounter{algorithmctr}[section]
\renewcommand{\thealgorithmctr}{\thesection.\arabic{algorithmctr}}
\newenvironment{algdesc}{\refstepcounter{algorithmctr}\begin{list}{}{%
			\setlength{\rightmargin}{0\linewidth}%
			\setlength{\leftmargin}{.05\linewidth}}%
		\rmfamily\small
		\item[]{\setlength{\parskip}{0ex}\hrulefill\par%
			\nopagebreak{\bfseries\textsf{Algorithm \thealgorithmctr~}}}}%
	{{\setlength{\parskip}{-1ex}\nopagebreak\par\hrulefill} \end{list}}
\journalname{Optimization and Engineering}
\begin{document}

\title{Minimizing Oracle-Structured Composite Functions
	\thanks{
		%Grants or other notes
		The work was supported in part by the National Key R\&D Program of China 
		with grant No. 2018YFB1800800, 
		by the Key Area R\&D Program of Guangdong Province 
		with grant No. 2018B030338001, 
		by Shenzhen Outstanding Talents Training Fund, 
		and by Guangdong Research Project No. 2017ZT07X152.
		Stephen Boyd’s work was funded in part by 
		the AI Chip Center for Emerging Smart Systems (ACCESS).
	}
}


\author{Xinyue Shen \and Alnur Ali \and Stephen Boyd}


\institute{Xinyue Shen \at
	Department of Electrical Engineering, Stanford University \\
	SRIBD and FNii,
	the Chinese University of Hong Kong \\
	\email{xinyues@stanford.edu}  
	\and
	Alnur Ali \at
	Department of Electrical Engineering, Stanford University \\
	Department of Statistics, Stanford University \\
	\email{alnurali@stanford.edu}
	\and
	Stephen Boyd  \at
	Department of Electrical Engineering, Stanford University \\
	\email{boyd@stanford.edu}
}

\date{Received: date / Accepted: date}

\maketitle

\begin{abstract}
We consider the problem of minimizing a composite convex function
with two different access methods: an \emph{oracle}, for which we can 
evaluate the value and gradient, and a \emph{structured function}, 
which we access only by solving a convex optimization problem.
We are motivated by two associated technological developments.  For the 
oracle, systems like PyTorch or TensorFlow can
automatically and efficiently compute gradients, given a computation
graph description.  For the structured function, systems like CVXPY
accept a high level domain specific language description of the problem,
and automatically translate it to a standard form for efficient solution.
We develop a method that makes minimal assumptions about the two functions,
does not require the tuning of algorithm parameters, and works well
in practice across a variety of problems.
Our algorithm combines a number of 
well-known ideas, including a low-rank quasi-Newton approximation of
curvature, piecewise affine lower bounds from bundle-type methods,
and two types of damping to ensure stability.
We illustrate the method on stochastic optimization, utility
maximization, and risk-averse programming problems, 
showing that our method is more efficient
than standard solvers when the oracle function contains
much data. 
\keywords{Composite convex optimization \and First-order oracles
	\and Structured optimization \and Quasi-second-order methods 
	\and Tuning-free methods}
\end{abstract}

\section{Introduction}\label{s-introduction}
Our story starts with a well studied problem, minimizing
a convex function that is the sum of two convex functions
with different access methods, referred to as a \emph{composite function}
\cite{nesterov2013gradient}.
The first function is smooth, and we can access it only by a few methods, such
as evaluating its value and gradient at a given point.
The other function is not necessarily smooth, but is structured, and 
we can access it only by solving a convex optimization problem that involves it.
In the typical setting, the second function is one for which we can 
efficiently compute its proximal operator, usually analytically.
Here we assume a bit more for the second function, specifically, 
that we can minimize it plus a structured function of modest complexity.

Our goal is to minimize such functions automatically, with a method 
that works well across a large variety of problem instances using its
default parameters.
We leverage new technological developments: systems for 
automatic differentiation that automatically and effficiently compute
gradients given a computation graph description, and
systems for solving convex optimization problems
described in a domain specific language for multiple parameter values.

\subsection{Oracle-structured composite function} \label{s-assumptions}
We seek to minimize $h(x) = f(x)+g(x)$ over $x\in \reals^n$.
We assume that
\BIT
\item $f:\Omega\to\reals$ is convex and differentiable, where 
$\Omega \subseteq \reals^n$ is convex and open.
We assume that a point $x^0 \in \Omega$ is known.
\item $g: \reals^n \to \reals\cup \{+\infty\}$ is convex, with closed sublevel
sets, and not necessarily differentiable.  Infinite values of $g$ encode
constraints on $x$.
\EIT
Our method will access these two functions in very specific ways.
\BIT
\item We can evaluate $f(x)$ and $\nabla f(x)$ at any $x$.
For $x\not\in \Omega$, our oracle returns the
value $+\infty$ for $f(x)$.
(We will discuss a few other possible access methods for $f$ in the sequel.)
\item We can minimize $g(x)$ plus another structured function.
Here we use the term structured loosely, to mean in the sense of Nesterov,
or as a practical matter, in a disciplined convex programming (DCP)
description.
This is an extension of the usual assumption in composite function
minimization that the proximal 
operator of $g$ can be evaluated analytically.
\EIT
We refer to $f$ as the oracle part of the objective,
$g$ as the structured part of the objective,
and $h=f+g$ as an oracle-structured composite objective.
We denote the optimal value of the problem as
\[
h^\star = \inf_x \left( f(x) + g(x) \right).
\]

\paragraph{Assumptions.}
We assume that the sublevel sets of $h$ are compact, 
and $h^\star < \infty$ (so at least some sublevel sets are nonempty).
These assumptions imply that $f(x)+g(x)$ has a minimizer, \ie, 
a point $x^\star$ with $h(x^\star)=h^\star$.
Our convergence proofs make the typical assumption that $\nabla f$ is Lipschitz 
continuous with constant $L$, but we stress that this is not used in the 
algorithm itself.

\paragraph{Optimality condition.}
The optimality condition is \BEQ
\label{e-opt-cond}
\nabla f(x)+ q = 0, \qquad q \in\partial g(x),
\EEQ
where $\partial g(x)$ denotes the subdifferential of $g$ at $x$.
For $x \in \Omega$ and $q \in \partial g(x)$, 
we can interpret $\nabla f(x)+q$ as a
residual in the optimality condition.
(We will use this in the stopping criterion of our algorithm.)
Our access to $f$ directly gives us $\nabla f(x)$; we will see that our
access method to $g$ indirectly produces a subgradient $q \in \partial g(x)$.

\subsection{Practical considerations}
We are motivated by two technological considerations related to 
our access methods to $f$ and $g$, which we mention briefly here.

\paragraph{Handling $f$.}
To handle $f$ we can rely on automatic differentiation systems
that have been developed in recent years, such as PyTorch \cite{NEURIPS2019_9015},
TensorFlow \cite{tensorflow2015-whitepaper}, and Zygote/Flux/Autograd 
\cite{maclaurin2015autograd,Flux.jl-2018,innes2019differentiable}.
Automatic differentiation is an old topic \cite{nolan1953analytical, adamson1969slang}
(see \cite{baydin2018automatic} for a recent review), but these
recent systems go way beyond the basic algorithms for automatic
differentiation, in terms of ease of use and run-time efficiency,
across multiple computation platforms.
We describe $f$ (but not its
gradient) using existing libraries and languages;
thereafter, $f(x)$ and $\nabla f(x)$ can be 
evaluated very efficiently, on many computation platforms, 
ranging from single CPU to multiple GPUs.
These systems are widely used throughout machine learning,
mostly for fitting deep neural networks.

\paragraph{Handling $g$.}
To handle $g$ we make use of domain specific languages (DSLs) 
for convex optimization, such as CVX~\cite{grant2014cvx}, 
CVXPY~\cite{cvxpy_paper,agrawal2018rewriting}, 
and Convex.jl~\cite{cvxjl}. These systems take a description of $g$
in a special language based on disciplined convex programming (DCP)~\cite{GBY:06}.
This description of $g$ is then automatically transformed into a 
standard form, such as a cone program, and then solved.
Recently, such systems have been enhanced to include parameters,
which are constants in the problem each time it is solved, 
but can be changed and the problem re-solved efficiently, skipping
the compilation process.
These systems are reasonably 
good at preserving structure in the problem during compilation,
so the solve times can be quite small when exploitable structure 
is present, which we will see is the case in our method.

We mention that $g$ can contain hidden additional variables.
By this we mean that $g$ has the form $g(x) = \inf_z G(x,z)$,
where $G$ is convex in $(x,z)$.  Roughly speaking, 
$z$ is the hidden variable that does not appear in $f$.  
Such functions
are immediately handled by structured systems, without 
any additional effort.   In particular, we do not need to work out
an analytical form for $g(x)$.
In this case just evaluating $g$ requires solving an optimization problem
(over $z$).  Our method will avoid any evaluations of $g$.

\paragraph{When to just use a structured solver.}
Finally, we mention that if $f$ is simple enough to be handled by
a DCP-based system, then simply minimizing $f+g$ using such a system
is the preferred method of solution.
We are interested here in problems where this is not the case.
Typically this means that $f$ is complex in the sense of 
involving substantial data, for example, a sample average of some
function with $10^6$ or more samples.
(We will see this phenomenon in the numerical examples given in 
\S\ref{s-experiments}.)

\paragraph{Contribution.}
The method we propose in the next section is in the family of
variable metric bundle methods, and closely related to a number of 
other methods found in the literature (we review related 
work in \S\ref{s-related-work}).  As an algorithm, our method is not 
particularly novel; we consider our contribution to be its careful 
design to be compatible with our access methods,
and the efficiency of the
method compared with standard solvers
when the function $f$ contains a large amount of data.

\section{Oracle-structured minimization method} \label{s-method}
In this section we propose a generic method for solving the oracle-structured
minimization problem,
which we call \emph{oracle-structured minimization method} (OSMM).
OSMM combines several well known methods from optimization,
including variable metric or quasi-Newton curvature estimates 
to accelerate convergence,
bundle methods that build up a piecewise affine model, and two types of 
damping, based on a trust penalty and a line search.
These are chosen to be compatible with our access methods.

We will denote the iterates with a superscript, so $x^k$ denotes the $k$th
iterate of the algorithm.
We will let $x^{k+1/2}$ denote the tentative iterate at the $(k+1)$st iteration,
before the line search.
We will assume that $x^0\in\Omega$, \ie, $f(x^0)<\infty$.  Our algorithm
will guarantee that $x^k \in \Omega$ for all $k$.
It is a descent method, \ie, $h(x^{k+1}) < h(x^k)$.  While $h(x^0)=\infty$
is possible, we will see that $h(x^k)< \infty$ for $k\geq 1$.

As with many other optimization algorithms, OSMM is 
based on forming an approximation of the function $f$ in each
iteration.

\subsection{Approximation of the oracle function}\label{s-oracle_approximation}
In iteration $k$, we form a convex approximation of $f$, given by 
$\hat f_k: \reals^n \to \reals \cup \{\infty\}$, based 
on information obtained from previous iterations and possibly prior 
knowledge of $f$.  
The approximation has the specific form
\BEQ\label{e-fhat}
\hat f_k(x) = l_k(x) + (1/2)(x-x^k)^T H_k (x-x^k).
\EEQ
Here $H_k$ is positive semidefinite,
and $l_k:\reals^n \to \reals \cup \{\infty\}$ is a convex
minorant of $f$, \ie, $l_k(x) \leq f(x)$ for all $x$.

\paragraph{Assumptions on $H_k$.}
The only assumption we make about $H_k$ is that it is positive
semidefinite and bounded, \ie, there exists a $C$ such that
$\|H_k \|_2 \leq C$ for all $k$.
In practice, $H_k$ accelerates convergence 
by serving as an estimate of the curvature of $f$.  
The simplest choice, $H_k=0$, results in 
an algorithm that converges
but does not offer the practical benefit of 
convergence acceleration.

\paragraph{Choice of $H_k$.}
There are many ways of choosing a positive semi-definite $H_k$ 
to approximate the curvature of $f$ at $x^k$.
One obvious choice is the Hessian $H_k = \nabla^2 f(x^k)$,
but this requires that $f$ be twice differentiable, and also 
violates our assumption about how we access $f$.
A lesser violation of the access method might use an
approximation of the Hessian based on evaluations of 
the mapping $z \mapsto \nabla^2f(x^k) z$ (\ie, Hessian-vector multiplication),
which can be practical in many cases
\cite{erdogdu2015convergence}.
A simple and effective choice is
$H_k = (a_k/n) I$, where $a_k$ is an approximation of $\Tr \nabla^2 f(x^k)$,
obtained for example by the Hutchinson 
method \cite{hutchinson1989stochastic,meyer2021hutch++}.

Quasi-Newton methods are a general class of curvature approximations 
that are compatible with our assumptions 
on the access method for $f$.
These methods, which have a very long history,
build up an approximation of $H_k$ using only the current 
and previously evaluated gradients
\cite{davidon1959variable,fletcher1963rapidly,
	broyden1965class,dennis1977quasi,byrd1994representations}.
When $H_k$ is low rank, or diagonal plus low rank, the method is 
practical even for large values of $n$.  (Such methods are often called
limited memory, since they do not require the storage of 
an $n \times n$ matrix.)
For OSMM we propose to use the
low-rank quasi-Newton choice given in \cite{fletcher2005new},
described in detail in \S\ref{appendix:s-curvature}.
We can express $H_k$ as
\BEQ\label{PSD_low_rank}
H_k = G_k G_k^T,
\EEQ
where $G_k\in\reals^{n\times r}$, and $r$ is a chosen (maximum) rank
for $H_k$.

As many others have observed, limited memory quasi-Newton methods
deliver most of their benefit for relatively small values of $r$, 
like $r=10$ or $r=20$.   These values allow the methods to be used 
even when $n$ is very large (say, $10^5$), since the storage requirement
(specifically, of $G_k$) grows linearly with $r$, and the computational
cost of evaluating $x^{k+1/2}$ grows quadratically in $r$,
and only linearly in $n$.

\paragraph{Assumptions on $l_k$.}
We make the usual assumption on the minorant $l_k$ that it is tight
at $x^k$, \ie, $l_k(x^k)=f(x^k)$.
It follows that $\hat f_k(x^k) = f(x^k)$.
It also follows that $l_k$ is differentiable at $x^k$, and 
$\nabla l_k(x^k) = \nabla f(x^k)$.
To see this, we note that since $l_k$ is a minorant of $f$,
tight at $x^k$, we have 
\[
\partial l_k(x^k) \subseteq \partial f(x^k) = \{ \nabla f(x^k) \}.
\]
The first inclusion can be seen since any affine lower bound on $l_k$,
tight at $x^k$, is also an affine lower bound on $f$, tight at $x^k$,
so its linear part is a subgradient of $f$ at $x^k$.
The right-hand equality holds since $f$ is differentiable, so its 
subdifferential contains only one element, its gradient.
Finally, since $\partial l_k(x^k)$ contains only $\nabla f(x^k)$,
we conclude it is differentiable at $x^k$, with gradient $\nabla f(x^k)$.

In addition to the mathematical assumptions about $l_k$ described above,
we will assume that $l_k$ has a structured description.
This implies that $\hat f_k$ has a structured description.

\paragraph{Minorants.}
The simplest minorant is the first order Taylor approximation
\[
l_k(x) = f(x^k) + \nabla f(x^k)^T(x-x^k).
\]
A more complex minorant is the piecewise affine minorant
\[
l_k(x) = \max_{i=1,\ldots, k} \left( f(x^i) + \nabla f(x^i)^T(x-x^i)\right),
\]
which uses all previously evaluated gradients of $f$.

For OSMM we propose the piecewise affine minorant
\BEQ\label{l_k_PWA}
l_k(x) = \max_{i = \max\{0, k-M+1\},\ldots,k} \left(
f(x^i) + \nabla f(x^i)^T(x-x^i) \right),
\EEQ
the pointwise maximum of the affine minorants from the previous
$M$ gradient evaluations, where $M$ is the memory.
With memory $M=1$, this reduces to the Taylor approximation.

The problem of choosing the memory $M$ is very similar to the problem of
choosing $r$, the rank of the curvature approximation.  The storage 
requirements grows linearly with $M$, and the computational cost
of evaluating $x^{k+1/2}$ grows quadratically with $M$.
As with the choice of $r$, small values such as $M=10$ or $M=20$
seem to work well in practice.

We mention a few additional useful minorants that use additional 
prior information about $f$ or its domain $\Omega$.
First, we can add to $l_k$ constraints that contain $\Omega$.
Suppose we know that $\tilde \Omega \supset \Omega$, where 
$\tilde \Omega$ has a structured description, \eg, a box.
We can then use the minorant
\[
l_k(x) = \max_{i = \max\{0, k-M+1\},\ldots,k} \left(
f(x^i) + \nabla f(x^i)^T(x-x^i) \right)
+ I_{\tilde \Omega}(x),
\]
where $I_{\tilde \Omega}$ is the indicator function of $\tilde \Omega$.
In a similar way, if a (constant) lower bound $\ell$ on $f$ is known, 
we can replace any minorant $l_k$ with $\max\{l_k(x),\ell\}$.

If it is known that $f$ is $\rho$-strongly convex, 
we can strengthen the piecewise affine minorant \eqref{l_k_PWA} 
to the piecewise quadratic minorant
\[
l_k(x) = \max_{i=\max\{0, k-M+1\},\ldots,k} \left( f(x^i) 
+ \nabla f(x^i)^T(x-x^i) + (\rho/2)\|x-x^i\|_2^2 \right).
\]
(Each term in the maximum has the same quadratic part $(\rho/2)\|x\|_2^2$,
so $l_k$ can be expressed as a piecewise affine function plus 
$(\rho/2)\|x\|_2^2$.)

\paragraph{Lower bound.}
We observe that
\[
\ell_k = \inf_{x} \left( l_k(x) + g(x) \right)
\]
is a lower bound on the optimal value $h^\star$.
It can be computed using a system for structured optimization.
At iteration $k$ we let $L_k$ denote the best (largest) lower bound
found so far,
\BEQ\label{e-Lk}
L_k = \max\{\ell_1, \ldots, \ell_k\}.
\EEQ

\subsection{Tentative update}\label{s-tentative-update}
At iteration $k$,
our tentative next iterate $x^{k+1/2}$ is obtained by minimizing
our approximation of $f$, plus $g$ and a trust penalty term:
\BEQ\label{e-xk+half}
x^{k+1/2} = \argmin_{x} \left( 
\hat f_k(x) + g(x) + \frac{\lambda_k}{2} \|x-x^k\|_2^2
\right).
\EEQ
The last term is a (Levenberg-Marquardt or proximal) trust penalty, 
which penalizes deviation from $x^k$;
the positive parameter $\lambda_k$ scales the trust penalty.
We assume that $x^{k+1/2}$ in \eqref{e-xk+half} can be computed
using a system for structured optimization.
(The minimizer in \eqref{e-xk+half} exists and is unique,
so $x^{k+1/2}$ is well defined.
To see this we observe that the objective is finite for $x=x^k$, and the
function being minimized is strictly convex.)
We note that while $\hat f_k(x^{k+1/2})$ and $g(x^{k+1/2})$ are finite,
$f(x^{k+1/2})=\infty$ (and therefore also $h(x^{k+1/2})=\infty$) is possible.

The two quadratic terms in the objective in \eqref{e-xk+half}
can be combined to express the tentative update as
\BEQ\label{e-xk+half2}
x^{k+1/2} = \argmin_{x} \left( 
l_k(x) + g(x) + \frac{1}{2}(x-x^k)^T (H_k + \lambda_k I)(x-x^k)
\right),
\EEQ
which shows that the trust penalty term can be interpreted as a
regularizer for the curvature estimate $H_k$.

In the DCP description of the problem \eqref{e-xk+half2},
using \eqref{PSD_low_rank} we express the last term as
\[
\frac{1}{2}(x-x^k)^T (H_k + \lambda_k I)(x-x^k) = 
\frac{1}{2}\|G_k^T(x - x^k)\|_2^2 + \frac{\lambda_k}{2} \|x - x^k\|_2^2.
\]
This keeps the problem \eqref{e-xk+half2} tractable 
when $r \ll n$ and $n$ is large.
In particular, there is no need to form the $n \times n$ matrix $H_k$.

\paragraph{Tentative update optimality condition.}
For future reference, we note that the optimality 
condition for the minimization in \eqref{e-xk+half2} that 
defines $x^{k+1/2}$ is
\BEQ\label{e-xk+half-opt}
0 \in \partial l_k(x^{k+1/2}) + \partial g(x^{k+1/2})
+ (H_k+\lambda_k I) (x^{k+1/2} - x^k).
\EEQ
When $l_k$ is the piecewise affine minorant \eqref{l_k_PWA},
its subdifferential $\partial l_k(x^{k+1/2})$ has the form
\BEQ \label{e-subdiff-lk}
\partial l_k(x^{k+1/2}) = \Co \{ \nabla f(x^i) \mid l_k(x^{k+1/2}) =
f(x^i)+ \nabla f(x^i)^T(x^{k+1/2}-x^i) \},
\EEQ
the convex hull of the gradients associated with the 
active terms in maximum defining $l_k$.  
In the simplest case when $l_k$ is differentiable at 
$x^{k+1/2}$, \ie, only one term is active, this reduces to
$\{ \nabla l_k(x^{k+1/2})\} = 
\{ \nabla f(x^i)\}$, where $i$ is the (unique) index for which
$l_k(x^{k+1/2}) = f(x^i)+ \nabla f(x^i)^T(x^{k +1/2}-x^i)$.

Once $x^{k+1/2}$ is computed 
by solving problem \eqref{e-xk+half},
we can recover specific subgradients in
the subdifferentials $\partial l_k(x^{k+1/2})$ and $\partial g(x^{k+1/2})$
that satisfy \eqref{e-xk+half-opt}.
%(How to do this is explained in \S\ref{appendix-q_k}.
%We only need to find one of the two subgradients, since 
%the three terms in \eqref{e-xk+half-opt} sum to zero,
%and we know the third term.)
%We will denote the specific subgradient in $\partial g(x^{k+1/2})$
%in \eqref{e-xk+half-opt} as $q^{k+1}$.
As explained in \S\ref{appendix-q_k}, 
the subgradient in $\partial l_k(x^{k+1/2})$
has the form $\sum_i \gamma_i \nabla f(x^i)$,
where $\gamma_i $ are nonnegative and sum to one,
and positive only for $i$ associated with 
active terms in the maximum that defines $l_k(x^{k+1/2})$.
The specific subgradient in $\partial g(x^{k+1/2})$ is
\BEQ\label{e-qk}
q^{k+1} = 
- \sum_i \gamma_i \nabla f(x^i)
-(H_k+\lambda_k I) (x^{k+1/2}-x^k) \in \partial g(x^{k+1/2}),
\EEQ
which will be useful in a stopping criterion,
as we shall see later in \S\ref{s-stopping-crit}.

\subsection{Descent direction}\label{s-descent_direction}
If $x^k$ is a fixed point of the tentative update,
\ie, $x^{k+1/2} =x^k$, 
then $x^k$ is optimal. To see this,
if $x^{k+1/2} =x^k$, from \eqref{e-xk+half-opt} we have
\BEQ\label{e-xk+half-opt-cond}
0\in \partial l_k(x^{k}) +  \partial g(x^{k}) = \nabla f(x^k) +  \partial g(x^{k}),
\EEQ
so $x^k$ is optimal. 
From the first inclusion in \eqref{e-xk+half-opt-cond}, we can also
conclude that $x^k$ minimizes $l_k(x)+g(x)$, so
$L_k = l_k(x^k) + g(x^k) = f(x^k) + g(x^k)$,
\ie, the lower bound in \eqref{e-Lk} is tight when $x^{k+1/2} = x^k$.

If $x^{k+1/2} \neq x^k$, the tentative step
\BEQ\label{e-vk}
v^k = x^{k+1/2}-x^k
\EEQ 
is a descent direction for $h$ at $x^k$, \ie,
for small enough $t>0$ we have $h(x^k + t v^k)<h(x^k)$.
That is, the directional derivative $h'(x^k; v^k)$ is negative.

To see this, we first observe that by \eqref{e-xk+half2},
\BEQ\label{desc_dir_step1}
l_k(x^{k+1/2}) + g(x^{k+1/2}) + \frac{1}{2}(v^k)^T(H_k+\lambda_k I)v^k
<
l_k(x^k) + g(x^{k}),
\EEQ
since $x^{k+1/2}$ minimizes the left-hand side, and the right-hand side
is the same expression, evaluated at $x^k$.
We also have
\[
l_k(x^{k+1/2}) + g(x^{k+1/2}) \geq
l_k(x^{k}) + g(x^{k}) + (l_k+g)'(x^k; v^k).
\]
Combining these two inequalities we get
\[
(l_k+g)'(x^k; v^k) <
- \frac{1}{2}(v^k)^T(H_k+\lambda_k I)v^k.
\]
Finally, we observe that
\[
(l_k+g)'(x^k; v^k) = (f+g)'(x^k; v^k) = h'(x^k; v^k),
\]
since $l_k$ is differentiable at $x^k$, with $\nabla l_k(x^k) = \nabla f(x^k)$.
So we have
\BEQ\label{e-descent}
h'(x^k; v^k) < -\frac{1}{2} (v^k)^T (H_k+\lambda_k I)v^k ,
\EEQ
which shows that $v^k$ is a descent direction for $h$ at $x^k$.

\subsection{Line search}\label{s-line-search}
The next iterate $x^{k+1}$ is found as 
\BEQ\label{e-xk+1}
x^{k+1} = x^k + t_k v^k = x^k + t_k (x^{k+1/2}-x^k),
\EEQ
where $t_k \in (0,1]$ is the step size.
When $t_k=1$, we say the step is un-damped.
We will choose $t_k$ using a variation on 
a traditional Armijo-type line search \cite{armijo1966minimization} 
that avoids additional evaluations of $g$.

For $t \in [0,1]$ we define
\BEQ\label{e_phi_k}
\phi_k(t) = f(x^k + t v^k) + t g(x^{k+1/2}) + (1-t)g(x^k).
\EEQ
Since the second and third terms are the chord above $g$, 
we have, for $t \in [0,1]$,
\BEQ\label{e-chord}
\phi_k(t) \geq h(x^k+t v^k).
\EEQ
Evidently $\phi_k(0)=h(x^k)$, and $\phi_k$ is differentiable, with
\[
\phi_k'(0) = \nabla f(x^k)^T v^k + g(x^{k+1/2})-g(x^k).
\]
Since $\nabla f(x^k) = \nabla l_k(x^k)$ and $l_k$ is convex, we get
\[
\nabla f(x^k)^T v^k = \nabla l_k(x^k)^T v^k \leq l_k(x^{k+1/2}) - l_k(x^k),
\]
so
we have
\[
\phi_k^\prime(0) \leq l_k(x^{k+1/2}) - l_k(x^k)+ g(x^{k+1/2})-g(x^k).
\]
Combining this with \eqref{desc_dir_step1}, we obtain
\BEQ\label{phi_prime_bound}
\phi_k^\prime(0) < - \frac{1}{2}(v^k)^T ( H_k + \lambda_kI) v^k.
\EEQ
Thus for $t>0$ small,
\BEQ\label{e-phik-bnd}
\phi_k(t) = \phi_k(0) + \phi_k(0)^\prime t + o(t^2) 
< h(x^k)  -\frac{t}{2}(v^k)^T ( H_k + \lambda_kI) v^k + o(t^2).
\EEQ

\paragraph{Step length.}
Let $\alpha, \beta \in (0,1)$.  We take $t_k = \beta^{j}$, 
where $j$ is the smallest nonnegative integer for which
\BEQ\label{armijo_step}
\phi_k(t_k) \leq h(x^k)
- \frac{\alpha t_k}{2} (v^k)^T (H_k+\lambda_k I)v^k
\EEQ
holds.  
(The condition~\eqref{armijo_step} holds for some $j$ by \eqref{e-phik-bnd}.)
A nice feature of this line search is that it does not require any additional
evaluations of the function $g$ (which can be expensive),
since we already know $g(x^k)$ and $g(x^{k+1/2})$.
As has been noted by many authors, the choice of the line 
search parameters $\alpha$ and $\beta$ is not critical.
Traditional default values such as 
\BEQ\label{e-line-search-params}
\alpha=0.05, \qquad \beta=0.5
\EEQ
work well in practice.

\subsection{Adjusting the trust parameter}
We have already observed that $\lambda_{k}$ is a regularizer for $H_{k}$.
A natural choice is to choose the regularizer parameter roughly 
proportional to $\tau_k= \Tr H_{k}/n$,
since $\tau_k I$ is the minimum Frobenius norm approximation 
of $H_{k}$ by a multiple of the identity. 
Thus we take
\BEQ\label{e-lamk}
\lambda_{k} = \mu_k \left(\tau_k + \tau_{\min} \right),
\EEQ
where $\mu_k$ gives the trust parameter relative to $\tau_k$,
and $\tau_{\min}$ is a positive lower limit.

We update $\mu_k$ by decreasing it when the line search is undamped,
\ie, $t_k=1$, and increasing it when the line search is damped,
\ie, $t_k<1$.
We do this with
\BEQ\label{e-muk}
\mu_{k+1} = \left\{
\begin{array}{ll} 
	\max\{ \gamma_\text{dec} \mu_k, \mu_{\min} \} & t_k=1 \\
	\min\{ \gamma_\text{inc}\mu_k, \mu_{\max} \} & t_k<1,
\end{array} \right.
\EEQ
where $\mu_\text{min}$ and $\mu_\text{max}$ are positive lower
and upper limits for $\mu_k$,
$\gamma_\text{dec} \in (0,1)$ is the factor by which we decrease
$\mu_k$,
and $\gamma_\text{inc} \in (1,\infty)$ is the factor by which we increase
$\mu_k$.
The values
\BEQ\label{e-lambda-update-params}
\tau_{\min} = 10^{-3}, \quad
\gamma_\text{dec} = 0.8, \quad 
\gamma_\text{inc} = 1.1, \quad
\mu_{\min} = 10^{-4}, \quad
\mu_{\max} = 10^5
\EEQ
give good results for a wide range of problems.
We can take $\mu_0=1$.

We mention one initialization that is useful when
$f$ is twice differentiable and we have the ability to evaluate
$z \mapsto \nabla^2f(x^k) z$ (\ie, Hessian-vector multiplication).
In this case we can replace $\tau_0$ with an estimate of 
$\Tr \nabla^2 f(x^0)/n$ obtained using the Hutchinson method \cite{hutchinson1989stochastic}.

\subsection{Stopping criteria}\label{s-stopping-crit}
We use two stopping criteria, one based on a gap 
between upper and lower bounds on the optimal value,
and the other based on an optimality condition residual.
The gap condition is simple:
\BEQ\label{e-gap-stop}
h(x^k) - L_k \leq \epsilon_\mathrm{abs}^\mathrm{gap} 
+ \epsilon_\mathrm{rel}^\mathrm{gap}|h(x^k)|,
\EEQ
where $\epsilon_\mathrm{abs}^\mathrm{gap}$ and 
$\epsilon_\mathrm{rel}^\mathrm{gap}$ are positive
absolute and relative gap tolerances, respectively.
Evaluating $L_k$ can be almost as expensive as evaluating $x^{k+1/2}$,
but it is used only in the stopping criterion.  To reduce
this overhead, we evaluate $L_k$ only every ten iterations.
Reasonable values for the gap tolerances are
$\epsilon_\mathrm{abs}^\mathrm{gap} = 10^{-4}$ and
$\epsilon_\mathrm{rel}^\mathrm{gap} = 10^{-3}$.

The residual based stopping criterion is tested whenever we
take an undamped step, \ie, $t_k=1$.  In this case $x^{k+1}=x^{k+1/2}$,
and we obtain $q^{k+1} \in \partial g(x^{k+1})$ in \eqref{e-qk},
so $\nabla f(x^{k+1}) + q^{k+1}$ is a residual for the optimality condition
\eqref{e-opt-cond}.
The stopping criterion is
\BEQ\label{stopping_criterion_3}
\frac{1}{\sqrt{n}}\|\nabla f(x^{k+1}) + q^{k+1}\|_2 \leq
\epsilon_\mathrm{abs}^\mathrm{res} 
+ \epsilon_\mathrm{rel}^\mathrm{res} 
\left(\frac{1}{\sqrt{n}}\|\nabla f(x^{k+1})\|_2 + \frac{1}{\sqrt{n}}\|q^{k+1}\|_2 \right),
\EEQ
where 
$\epsilon_\mathrm{abs}^\mathrm{res}$ and 
$\epsilon_\mathrm{rel}^\mathrm{res}$ are
relative and absolute residual tolerances.
(Dividing the norm expressions above by $\sqrt n$ gives the
root mean square or RMS values of the argument.)
Reasonable values for these parameters are
$\epsilon_\mathrm{abs}^\mathrm{res}=10^{-4}$ and
$\epsilon_\mathrm{rel}^\mathrm{res}=10^{-3}$.

\subsection{Algorithm summary}
We summarize OSMM in algorithm \ref{alg:1}.
\begin{algdesc}\label{alg:1}
	\emph{Oracle-structured minimization method.}
	\begin{tabbing}
		{\bf given} an initial point $x^0\in\Omega$.\\
		{\bf for} $k=0, 1, \ldots, k_{\max}$\\
		\qquad \= 1.\ \emph{Form a surrogate objective.}
		Form $l_k$ and $H_k$.\\
		\> 2.\ \emph{Tentative step.} Compute $x^{k+1/2}$ by \eqref{e-xk+half}.\\
		\> 3.\ \emph{Line search and update.} Set line search 
		step size $t_k$ by \eqref{armijo_step}
		and $x^{k+1}$ by \eqref{e-xk+1}.\\
		\> 4.\ \emph{Compute lower bound.}  If $k$ is a multiple of $10$,
		evaluate $L_k$.\\
		\> 5.\ \emph{Check stopping criterion.} Quit if \eqref{e-gap-stop} or 
		\eqref{stopping_criterion_3} holds.\\
		\> 6.\ \emph{Update trust penalty parameter.} Update $\lambda_{k+1}$ by \eqref{e-lamk}
		and \eqref{e-muk}.
	\end{tabbing}
\end{algdesc}

The algorithm parameters in OSMM are the memory of the minorant $M$, the rank of 
the curvature estimate $r$, the line search parameters given in
\eqref{e-line-search-params}, the $\lambda_k$ update parameters
given on \eqref{e-lambda-update-params}, and the relative and absolute
gap and residual tolerances, given in \S\ref{s-stopping-crit}.
The practical performance of OSMM is not particularly sensitive to 
the choice of these parameters; our implementation uses as default values
the ones described above, with memory $M=20$ and rank $r=20$.

\subsection{Implementation}
We have implemented OSMM in an open-source Python package, available at
\begin{center}
	\url{https://github.com/cvxgrp/osmm}.
\end{center}
The user supplies a PyTorch description of the oracle function $f$, 
a CVXPY description of the structured function
$g$,  and an initial point $x^0\in\Omega$. 
The package invokes PyTorch to evaluate $f$ and its gradient $\nabla f$.  
The convex model $\hat f_k$ is then formed and handed off to CVXPY to 
efficiently compute the next tentative iterate.

\subsection{Related work}\label{s-related-work} There is a lot of prior work
related to the method proposed in this paper.  The work on variable metric
bundle
methods~\cite{lemarechal1978nonsmooth,fukushima1984descent,kiwiel1990proximity,schramm1992version,lemarechal1995new,mifflin1996quasi,lemarechal1997variable,mifflin1998quasi,lukvsan1998bundle,kiwiel2000efficiency,teo2010bundle,yu2010quasi,noll2013bundle,Bagirov2014,van2016strongly,de2016doubly,frangioni2002generalized,de2014convex,van2018incremental}
and (inexact) proximal Newton-type
methods~\cite{levitin1966constrained,schmidt2009optimizing,sra2012optimization,becker2012quasi,lee2014proximal,scheinberg2016practical,li2017inexact,ghanbari2018proximal,yue2019family,becker2019quasi,lee2019inexact,asi2019modeling,mordukhovich2020globally}
are probably the most closely related, though our method differs
in key ways arising from our assumed access methods.  Proximal Newton methods
are similar in philosophy to our approach, as both types of methods really shine
when the structured part of the objective can be minimized efficiently.
However, on a purely technical level, proximal Newton methods generally form
the (usual) first-order model of the oracle part of the objective, without 
a bundle term, 
which is of course different than the model we use,
and can affect the convergence speed.

A few variable metric bundle methods have been proposed over the years, but
these methods also differ from ours in subtle (but important) ways.  In
general, these methods are not geared towards solving composite minimization
problems.  
Additionally, our line search method is different, 
as ours is carefully designed to satisfy the access conditions.

Finally, there has been a surge of interest in scalable quasi-Newton
methods~\cite{erdogdu2015convergence,pmlr-v70-ye17a,gower2018tracking,huang2020span} recently.
The focus has been on developing low-rank, regularized, and sub-sampled
Newton-type methods.  Of these, the family of low-rank approximations to the
Hessian are most related to the curvature estimate used in our method.

\section{Convergence}
In this section we show the convergence of OSMM,
and give a numerical example to illustrate its typical performance,
as the two most important algorithm parameters, memory $M$ and rank $r$, vary.
(More examples will be given in \S\ref{s-experiments}.)

\subsection{Convergence}\label{s-convergence}
We demonstrate two types of convergence.  First, we show the iterates $x^k$ converge asymptotically in a certain sense, \ie, the sequence of tentative steps $v^k \to 0$ as $k \to \infty$.  
%(From \eqref{e-xk+half-opt}, this implies that when the iterates $x^k$ do converge, the convergence point must be optimal.)  
Second, we show the objective values generated by OSMM converge to the optimal value.  
Thus, any cluster point of the sequence of iterates $x^k$ is an optimal point.
The results additionally require the (standard) assumption that $\nabla f$ be Lipschitz continuous.

\paragraph{Convergence of iterates.}
To see that $v^k \to 0$, we observe from \eqref{e-chord}~and \eqref{armijo_step} that
\BEQ\label{e-quad-form-bound}
h(x^{k+1}) - h(x^k) \leq - \frac{\alpha t_k}{2} (v^k)^T (H_k + \lambda_k I) v^k \leq 0.
\EEQ
Since $h(x^k)$ is decreasing and bounded below by $h^\star$, it converges,
which implies that 
the left-hand side of \eqref{e-quad-form-bound} converges to zero as
$k \to \infty$.  This in turn implies that
\[
t_k (v^k)^T (H_k + \lambda_k I) v^k \to 0 \quad \mathrm{as} \quad k \to \infty.
\]
By construction $\lambda_k$ is lower bounded by $\tau_{\min} \mu_{\min}$,
and $t_k$ is also bounded away from $0$, 
as will be proved later in \S \ref{s-tk-one-io},
so we get that $v^k\to 0$, as claimed.

\paragraph{Convergence of objective values.}
We require the basic fact that there exists a subsequence of iterates $(k_\ell)$ accepting unit step length, \ie, that $t_{k_\ell} = 1$ for all $\ell$.  We give a proof of this fact in \S\ref{s-tk-one-io}, assuming that $\nabla f$ is Lipschitz continuous.

Now, by convexity,
\[
h(x^{k_\ell+1}) - h^\star \leq (\nabla f(x^{k_\ell+1}) + q^{k_\ell+1})^T (x^{k_\ell+1} - x^\star),
\]
for some $q^{k_\ell} \in \partial g(x^{k_\ell})$.  Adding and subtracting any $l_{k_\ell}^{k_\ell+1} \in \partial l_{k}(x^{k+1})$, we get\label{key}
\BEQ\label{e-func-val-descent0}
h(x^{k_\ell+1}) - h^\star \leq \big( ( \nabla f(x^{k_\ell+1}) - l_{k_\ell}^{k_\ell+1} ) + ( l_{k_\ell}^{k_\ell+1} + q^{k_\ell+1} ) \big)^T (x^{k_\ell+1} - x^\star).
\EEQ

From our assumptions on $h$, and because OSMM is a descent method, we see that $\| x^{k_\ell+1} - x^\star \|_2$ is bounded.  Moreover, because $H_k$ and $\lambda_k$ are uniformly bounded, we get from \eqref{e-xk+half-opt} that
\[
\| l_{k_\ell}^{k_\ell+1} + q^{k_\ell+1} \|_2 \leq \| H_{k_\ell} 
+ \lambda_{k_\ell} I \|_2 \| v^{k_\ell} \|_2 \lesssim \| v^{k_\ell} \|_2,
\]
where we write $a \lesssim b$ to mean that $a$ and $b$ satisfy $a \leq C b$, for some $C > 0$.  Finally, it can be shown (see \S\ref{s-minorant-accuracy} for details) that the limited memory piecewise affine minorant satisfies
\BEQ\label{e-minorant-accuracy}
\| \nabla f(x^{k_\ell+1}) - l_{k_\ell}^{k_\ell+1} \|_2 \lesssim \max_{j=\max\{0,\ldots,k_\ell-M+1\},\ldots,k_\ell} \| v^j \|_2.
\EEQ
Putting these together with \eqref{e-func-val-descent0}, we obtain that
\[
h(x^{k_\ell+1}) - h^\star \lesssim \max_{j=\max\{0,\ldots,k_\ell-M+1\},\ldots,k} \| v^j \|_2.
\]
Earlier we showed that $v^k \to 0$, so $h(x^{k_\ell+1}) \to h^\star$ as $\ell \to \infty$.  Because the sequence of objective values $h(x^k)$ is convergent, we get $h(x^k) \to h^\star$, as claimed.

\subsection{A numerical example}\label{s-numerical-cvg-example}
Next we investigate the convergence of OSMM numerically.  
Here and throughout this paper, 
we use a Tesla V100-SXM2-32GB-LS GPU with 32 gigabytes of memory 
to evaluate $f$ and $\nabla f$ via PyTorch, 
and an Intel Xeon E5-2698 v4 2.20 GHz CPU 
to compute the tentative iterates via CVXPY and for any baselines.

\paragraph{Problem formulation.}  
We consider an instance of the Kelly gambling problem~\cite[\S 4]{BoydVand04}
\cite{kelly2011new,busseti2016risk},
\BEQ\label{e-kelly}
\begin{array}{ll}
	\mbox{maximize} & \sum_{i=1}^N \pi_i \log ( r_i^T x ) \\
	\mbox{subject to} & x \geq 0, \quad \ones^T x = 1,
\end{array}
\EEQ
where $x \in \reals^n$ is the variable, and $r_i \in \reals^n_+$,
$i=1,\ldots,N$, $\pi \in \reals^N_+$ are the problem data, with $\sum_{i=1}^N \pi_i =1$.
%Here $x_j$ is the fraction of our wealth we put on bet $j$,
%$(r_i)_j$ gives the return of bet $j$ with outcome $i$, and 
%$\pi_i$ is the probability of outcome $i$.
Details about the Kelly gambling problem
	and what the variable and data represent
	can be found in~\S\ref{subsec:experiment-kelly}.
To put problem \eqref{e-kelly} into our oracle-structured form, 
we take $f$ to be the objective in~\eqref{e-kelly}, 
and $g$ to be the indicator function of the constraints,
in this case the probability simplex.

\paragraph{Problem instance.}
Our problem instance has $n=100$ bets and $N = 1{,}000{,}000$ 
possible outcomes.
%Note that simply storing the returns $r_i$,
%which are dense in this instance, requires roughly four gigabytes, 
%causing most existing solvers to struggle, as we will soon see.  
How the data were generated
	can also be found in~\S\ref{subsec:experiment-kelly}.
This problem instance requires 400 seconds to solve using 
MOSEK~\cite{mosek}, a high performance commercial solver, with accuracy set to
its lower value. 

\begin{figure}
	\centering
	\includegraphics[width=2.2 in]{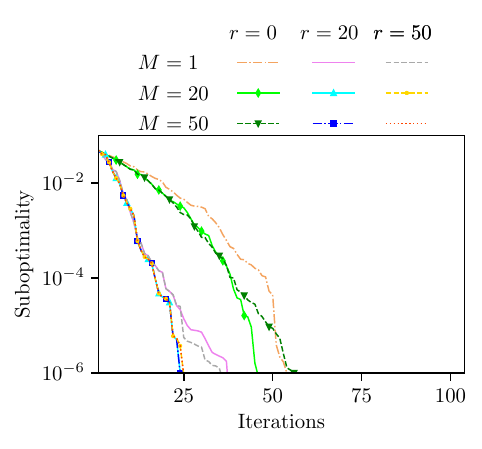}
	\includegraphics[width=2.07 in]{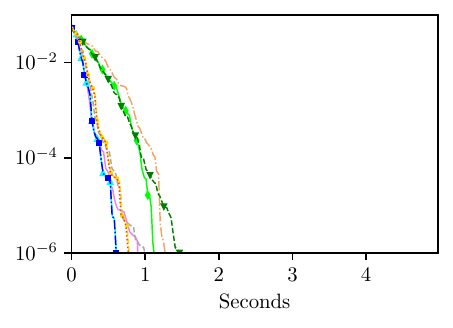}
	\includegraphics[width=2.2 in]{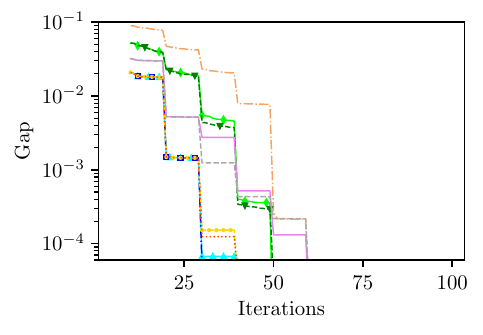}
	\includegraphics[width=2.07 in]{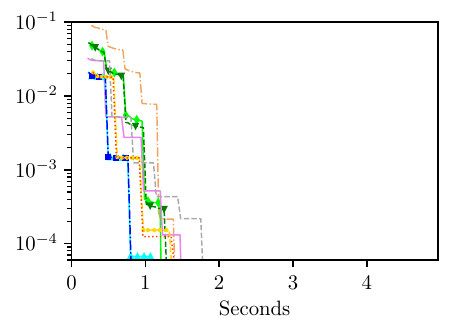}
	\includegraphics[width=2.2 in]{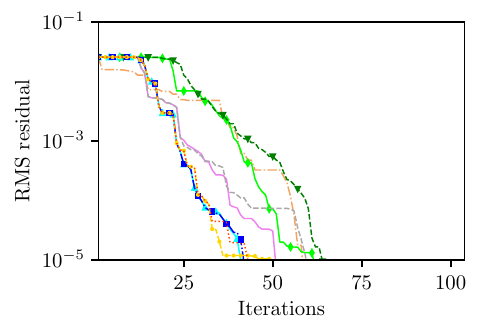}
	\includegraphics[width=2.07 in]{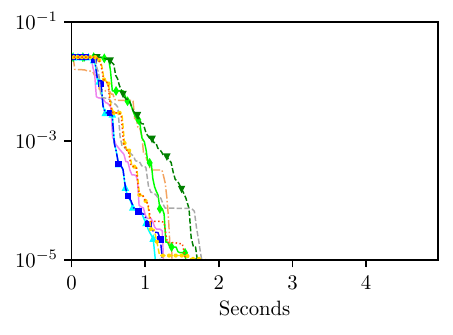}
	\caption{Suboptimality (top row), gap (middle row),
		and RMS residual (bottom row), versus
		iterations (left column) and run-time in seconds (right column) for
		OSMM on the Kelly gambling problem for the nine combinations of 
		rank $r$ and memory $M$. }
	\label{fig:convergence_example}
\end{figure}

\paragraph{Results.}  We solve~\eqref{e-kelly} by OSMM with rank $r \in
\{0,20,50\}$ and memory $M \in \{1,20,50\}$, a total of nine
choices of algorithm parameters.
The choices $r=0$ and $M=1$ correspond to using
no estimate of curvature, and no memory, respectively.
OSMM is run for 100 iterations for each combination of $r$ and $M$; $x^{k+1/2}$ and $L_k$ are computed using ECOS~\cite{Domahidi2013ecos}.
%The optimal objective value is -0.079.
Figure~\ref{fig:convergence_example}
shows the convergence of OSMM,
in terms of iterations (the left column)
and elapsed time (the right column).
The top row shows the true suboptimality $h(x^k) - h^\star$, which of course
we do not know when the algorithm is running, since we do not know $h^\star$.
The middle row shows the gap $h(x^k)- L_k$ 
(which we do know as the algorithm runs).   
The bottom row shows the RMS value of 
the optimality condition residual, \ie, the left-hand side
of~\eqref{stopping_criterion_3}.

The timing results averaged from $10$
	repetitions are shown in Table~\ref{tab:example_time_details}.
%More detail is given in table~\ref{tab:example_time_details}.
The total solve time and iterations are based on achieving a high accuracy
of $10^{-6}$, which is far more accurate than would be needed in practice.
For our nine choices of algorithm parameters, the total OSMM time 
ranges from around 0.87 to 1.4 seconds, substantially faster than using MOSEK
to solve the problem.

From the results in the figure and the table, 
as $r$ increases the convergence
becomes faster in terms of iterations, but with a larger 
$r$ the update of $x^{k+1/2}$ takes more time,
so each iteration becomes more expensive.
A good trade-off value is $r=20$.
The memory $M$ taking values $20$ and $50$
yields faster convergence in iterations and smaller gaps
than $M$ taking value $0$, but as $M$ increases
each iteration also becomes more expensive,
and the run-time efficiency is the best when $M=20$.
In summary, we see that rank value $r=20$ 
and memory value $M=20$ yield the fastest convergence
in terms of run-time.  We have observed this in other problem
instances as well, so these are the default values in OSMM.

\begin{table}
	\caption{Times to evaluate tentative update $x^{k+1/2}$ and lower bound $L_k$,
		total solve time and iterations to reach $10^{-6}$ accuracy, and average number of evaluations of $f$
		per iteration, for OSMM with the nine combinations of rank $r$ and memory $M$.}
	\label{tab:example_time_details}
	\centering
	\begin{tabular}{ll lll ll}
			\hline\noalign{\smallskip}
			\multicolumn{2}{l}{Rank and memory} &
			\multicolumn{2}{l}{Compute times (sec.)}&
			\multicolumn{3}{l}{Total solve time and iterations} \\
			\cmidrule(lr){1-2}
			\cmidrule(lr){3-4}
			\cmidrule(lr){5-7}
			$r$ & $M$ & $x^{k+1/2}$ & $L_k$ & Time (sec.) & Iterations & $f$ evals./iter.\\
			\noalign{\smallskip}\hline\noalign{\smallskip}
			0 & 1    & 0.016 & 0.0069 & 1.4 & 62 & 4.1\\
			0 & 20 & 0.019 & 0.0093 & 1.2 & 49 & 3.5\\
			0 & 50 & 0.023 & 0.015 & 1.3  & 49 & 3.4\\
			\cmidrule(lr){1-2}
			\cmidrule(lr){3-5}
			\cmidrule(lr){6-7}
			20 & 1    & 0.019  & 0.0069 & 0.88 & 36 & 2.9\\
			20 & 20 & 0.023  & 0.0093 & 0.87 & 33 & 2.9\\
			20 & 50 & 0.029  & 0.014  & 0.90 & 33 & 2.9\\
			\cmidrule(lr){1-2}
			\cmidrule(lr){3-5}
			\cmidrule(lr){6-7}
			50 & 1    & 0.021  & 0.0076& 0.92 & 32& 2.5\\
			50 & 20 & 0.027  & 0.010 & 0.88  & 29  & 2.5\\
			50 & 50 & 0.035  & 0.015 & 0.89 & 29 & 2.5\\
			\noalign{\smallskip}\hline
		\end{tabular}
\end{table}

\section{Generic applications}\label{s-applications}
In this section we describe some generic applications that reduce to our 
specific oracle-structured composite function minimization problem.
We focus on the function $f$, which should be differentiable, but not so
simple that it can be handled directly in a structured optimization
system.  We will see that this generally occurs when there is a lot of 
data required to specify $f$.

\subsection{Stochastic programming} 
\paragraph{Sample average approximation for stochastic programming.}
We start with the problem of minimizing $\Expect F(x,\omega) + g(x)$,
where $F$ is convex in $x$ for each value of the random variable $\omega$.
We will approximate the first objective term using a sample average.
We generate $N$ independent samples $\omega_1, \ldots, \omega_N$,
and take
\BEQ\label{e-emp-mean}
f(x) = \frac{1}{N} \sum_{i=1}^N F(x,\omega_i).
\EEQ
As a variation we can use importance sampling to get a lower variance estimate 
of $\Expect F(x,\omega)$. To do this we generate the samples from a 
proposal distribution with density $q$, and form the sample average estimate
\BEQ\label{e-imp-sampling}
f(x) = \frac{1}{N} \sum_{i=1}^N \frac{p(\omega_i)}{q(\omega_i)}F(x,\omega_i),
\EEQ
where $p$ is the density of $\omega$.

In both cases we can take $N$ to be quite large, 
since we will only need to evaluate $f$
and its gradient, and current systems for this are very efficient.
For example, the gradients $\nabla F(x,\omega_i)$ can be computed in
parallel.  As a practical matter, this happens automatically, with no
or very little directive from the user who specifies $f$.

\paragraph{Validation and stopping criterion tolerance.}
The functions $f$ given in \eqref{e-emp-mean} and \eqref{e-imp-sampling}
are only approximations of the true objective $\Expect F(x,\omega)$,
though we hope they are good approximations when we take $N$ large,
as we can with OSMM.
To understand how accurate the sample average \eqref{e-emp-mean} is,
we generate another set of independent samples 
$\tilde \omega_1,\ldots, \tilde \omega_N$
and define the validation function as
\[
f^\mathrm{val}(x) = \frac{1}{N} \sum_{i=1}^N 
F(x, \tilde \omega_i)
\]
(and similarly if we use importance sampling).
The magnitude of the difference $|f^\mathrm{val}(x)-f(x)|$ gives us a 
rough idea of the accuracy in approximating $\Expect F(x,\omega)$.
(Better estimates of the accuracy can be obtained by repeating
this multiple times, but we are interested in only a crude estimate.)

Solving the oracle-structured problem to an accuracy substantially better than
the Monte Carlo sampling accuracy does not make sense in practice.
This justifies replacing the absolute gap tolerance
$\epsilon_\mathrm{abs}^\mathrm{gap}$ with the maximum of 
a fixed absolute tolerance and 
the sampling error $|f^\mathrm{val}(x^k) - f(x^k)|$.
(We can evaluate the sampling error whenever we evaluate the gap,
\ie, every 10 iterations.)
Roughly speaking, we stop when we know we have solved the problem
to an accuracy that is better than our approximation.

\subsection{Utility maximization}
An important special case of stochastic optimization is
utility maximization, where we seek to maximize
\BEQ\label{e-util-max}
\Expect U(H(x,\omega)) - g(x),
\EEQ
where $U:\reals \to \reals$ is a concave increasing utility function,
and $H$ is concave in $x$ for each $\omega$.
The first term, the expected utility, is concave.
This is equivalent to the stochastic programming problem of 
minimizing
\[
\Expect (-U(H(x,\omega)) + g(x),
\]
which is stochastic programming with $F= -U\circ H$, which is convex
in $x$.  We can replace the expectation with a sample average
using \eqref{e-emp-mean}, or an importance sampling sample average 
using \eqref{e-imp-sampling}.
Utility maximization is a common method for handling the variance or 
uncertainty in a stochastic objective; it introduces risk aversion.

\subsection{Conditional and entropic value-at-risk programming}\label{s-EVaR}
Another method to introduce risk aversion into a stochastic optimization
problem is to mimimize value-at-risk (VaR), or a specific quantile of
$F(x,\omega)$, where
$F$ is convex in $x$ for each value of the random variable $\omega$.
The value-at-risk is defined as
\[
\VaR(F(x,\omega);\eta) = \inf \{ \gamma \mid \Prob(F(x,\omega) \leq \gamma) \geq \eta \},
\]
where $\eta \in (0,1)$ is a given quantile.
With an additional structured convex function $g$ in the objective, we obtain
the VaR problem
\begin{align}\label{prob-var}
\begin{array}{ll}
\mbox{minimize} & \VaR(F(x,\omega); \eta)) + g(x),
\end{array}
\end{align}
which, roughly speaking, is the problem of minimizing the $\eta$-quantile of
the random variable $F(x,\omega)+g(x)$.
Aside from a few special cases, this problem is not convex.  
(Recent work on VaR and
methods for the closely related problem 
of handling probability constraints can be found in, \eg, 
\cite{danielsson2013fat,danielsson2008optimal}.)
We proceed by replacing the nonconvex VaR term with a convex upper bound on VaR,
such as conditional value-at-risk (CVaR)
\cite{CVAR_uryasev_rockafellar,rockafellar2002conditional} 
or entropic value-at-risk (EVaR) \cite{ahmadi2011information}.
Beyond resulting in tractable convex problems, these upper bounds possess a number 
of nice properties, such as being coherent risk measures;
%which VaR is not; 
see \cite{rockafellar2002conditional,ahmadi2011information} for a discussion.

CVaR is given by
\BEQ\label{e-cvar}
\CVaR(F(x,\omega);\eta) = \inf_{\alpha \in \mathbf{R}} \Bigg\{ \frac{\Expect (F(x,\omega) - \alpha)_+}{1-\eta} + \alpha \Bigg\},
\EEQ
where $(z)_+ = \max\{z,0\}$, and EVaR is given by
\BEQ\label{e-evar}
\EVaR(F(x,\omega);\eta) = \inf_{\alpha > 0} \Bigg\{ \alpha \log \Bigg( \frac{\Expect \exp(F(x,\omega)/\alpha)}{1-\eta} \Bigg) \Bigg\}.
\EEQ
Both of these are convex functions of $x$.
We have, for any $x$, 
\[
\VaR(F(x,\omega);\eta) \leq \CVaR(F(x,\omega);\eta) \leq \EVaR(F(x,\omega); \eta)
\]
(see, \eg, \cite{ahmadi2011information}).
With CVaR, we obtain the convex problem
\[
\begin{array}{ll}
\mbox{minimize} & \CVaR(F(x,\omega); \eta)) + g(x),
\end{array}
\]
and similarly for EVaR.

We now show how the CVaR and EVaR problems can be approximated as oracle-structured
problems.
We start with CVaR.
We generate independent samples $\omega_i$, $i=1, \ldots, N$,
and replace the expectation with the empirical mean,
\[
f(x,\alpha) = \frac{1}{N} \sum_{i=1}^N
\frac{(F(x,\omega_i) - \alpha)_+}{1-\eta} + \alpha,
\]
with variables $x$ and $\alpha$.   This function is jointly convex in $x$ and $\alpha$;
minimizing over $\alpha$ gives $\CVaR(F(x,\omega);\eta)$ for the empirical 
distribution.
We adjoin $\alpha$ to $x$ to obtain a problem in oracle-structured form, \ie,
minimizing $f(x,\alpha)+g(x)$, over $x$ and $\alpha$.
(That is, we take $(x,\alpha)$ as what we call $x$ in our general form.)
While $f$ is not differentiable in $(x,\alpha)$, we have observed that 
our method still works very well.

For EVaR, we obtain
\[
f(x,\alpha) = \alpha \log \Bigg( \frac{1}{N} \sum_{i=1}^N
\frac{\exp(F(x,\omega_i)/\alpha)}{1-\eta} \Bigg),
\]
which is jointly convex in $x$ and $\alpha$.
(To see convexity, we observe that $f$ is the perspective function
of the log-sum-exp function; see \cite[\S 3.2.6]{BoydVand04}.)
As with CVaR,
minimizing over $\alpha$ yields $\EVaR(F(x,\omega);\eta)$ for the empirical
distribution.
Unlike our approximation with CVaR, this function is differentiable.

\subsection{Generic exponential family density fitting} \label{s-exp-fam-fit}
We consider fitting an exponential family of densities, given by
\BEQ\label{e-exp-fam}
p_\theta(z) = e^{- \left(\phi(z)^T \theta + A(\theta) \right)}, %\quad z\in\mathcal{S}
\EEQ
to samples $z_1,\ldots, z_m\in\mathcal{S}$.  Here, $\mathcal S$ is the support of the density, $\phi : \mathcal{S} \to \reals^n$ is the sufficient statistic, and $\theta \in \Theta \subseteq \reals^n$ is the canonical parameter (to be fitted).  The density normalizes via the log-partition (or cumulant generating) function
\[
A (\theta) = \log \int_\mathcal{S} e^{-\phi(z)^T \theta} \; dz.
\] 
We assume $\Theta$ is convex, and additionally that it only contains parameters for which the log-partition function is finite.

The negative log-likelihood, given samples $z_1,\ldots,z_m$, is
\[
\sum_{i=1}^m \log p_\theta(z_i) = c^T \theta + m A(\theta),
\]
where $c = \sum_{i=1}^m \phi(z_i)$.  So maximum likelihood estimation of $\theta$ corresponds to solving the density fitting problem
\begin{align}\label{exp_density_prob_0}
\begin{array}{ll}
\mbox{minimize} & \frac{1}{m}c^T \theta + A(\theta) + \lambda r(\theta)\\
\mbox{subject to} 
& \theta\in\Theta,
\end{array}
\end{align}
with variable $\theta$.  We can include a (potentially nonsmooth) convex regularization term $\lambda r(\theta)$ in the objective, where $\lambda \geq 0$ is the regularization strength, and $r : \Theta \to \reals$ is the regularizer.  Since the log-partition function is convex~\cite{wainwright2008graphical}, the density fitting problem~\eqref{exp_density_prob_0} is also convex.

The log-partition function is generally intractable except for a few special cases, so we replace the integral in $A(\theta)$ with a finite sum using importance sampling, \ie,
\[
\int_\mathcal{S} e^{-\phi(z)^T \theta} \; dz \approx
\frac{1}{N} \sum_{i=1}^N \frac{1(\omega_i\in\mathcal{S})}{q(\omega_i)}e^{-\phi(\omega_i)^T \theta},
\]
where $\omega_i, \; i=1,\ldots,N$, are independent draws from the proposal distribution $q$.  
The number of samples $N$ can be very large, especially when the number of dimensions $n$ is moderate.  
When $\mathcal S$ is bounded and its dimension is small, we can simply use a Riemann sum, so that the samples $\omega_i$ are lattice points in $\mathcal S$ and we have $q(\omega_i) = 1/|\mathcal{S}|$.
The problem~\eqref{exp_density_prob_0} is clearly in oracle-structured form, once we take $A(\theta)$ (with its Monte Carlo approximation) to be $f$, and the rest of the objective plus the indicator of $\Theta$ to be $g$.

A number of interesting regularizers are possible.  The squared $\ell_2$ norm, \ie, $r(\theta) = \| \theta \|_2^2$, is of course a natural choice.  
When $\mathcal{S}$ is bounded, a different option is to use the squared $L_2$ norm of the gradient of the log-density
\[
r(\theta) = \int_\mathcal{S} \|\nabla \log p_\theta(z)\|_2^2 \; dz.
\]
This regularizer enforces a certain kind of smoothness: the regularized density $p_\theta$ tends to the uniform distribution on $\mathcal S$, as the regularization strength $\lambda$ grows.  Finally, observe that we can write
\[
\int_\mathcal{S} \|\nabla \log p_\theta(z)\|_2^2 \; dz =
\int_\mathcal{S} \|D\phi(z)^T \theta\|_2^2 \; dz
= \theta^T \left(\int_\mathcal{S} D\phi(z) D\phi(z)^T \; dz \right) \theta,
\]
where $D\phi$ is the Jacobian of sufficient statistic $\phi$.  We can replace the integral with a finite sum again, and obtain the regularizer
\begin{equation*}
r(\theta) \approx \theta^T \left(\frac{1}{N} \sum_{i=1}^N 
\frac{1(\omega_i\in\mathcal{S})}{q(\omega_i)} D\phi(\omega_i) D\phi(\omega_i)^T \right) \theta.
\end{equation*}
% \EEQ

\section{Numerical examples}\label{s-experiments} In this section we
demonstrate the performance of OSMM through several numerical examples, 
all taken from the
generic applications described in the previous section.
We start by showing results for two different portfolio selection problems,
the Kelly gambling example (shown earlier in \S\ref{s-numerical-cvg-example}),
and minimizing CVaR.  We then show
a density estimation example.  After that we present results for a
supply chain optimization problem with entropic value-at-risk minimization.  
All of these examples are structured stochastic optimization problems, and we
use simple sample averages to approximate expectations.
OSMM is designed to handle the case when $f$ is complex, 
which in the case of sample averages means $N$ is large.
We will see that when $N$ is small, OSMM is actually slower than just 
solving the problem directly using a structured solver; when $N$ is large,
it is much faster (and in many cases, directly using a structured solver
fails).

We report the time needed for OSMM to reach high accuracy, \ie, $h(x^k) - h^\star \leq 10^{-6}$.  
(This makes a fairer comparison with MOSEK, SCS~\cite{ocpb:16,scs}, and ECOS.)
We also indicate when practical accuracy is reached, using our default 
parameters, and the sampling accuracy.
We use the default parameters in OSMM,
and use ECOS as the solver in CVXPY to compute the tentative 
iterate $x^{k+1/2}$.  
We do not perform any parameter tuning for our method.

\subsection{Kelly gambling}\label{subsec:experiment-kelly}

\paragraph{Problem formulation.} In the Kelly gambling problem there are $n$
bets we can wager on, and $N$ possible outcomes, with probabilities
$\pi_i$, $i=1, \ldots, N$.
The bet returns are given by $r_i\in \reals_+^n$, where $(r_i)_j$ is
the return, \ie, the amount you win for each dollar you put on bet $j$
when outcome $i$ occurs.
We are to choose $x\in \reals^n$, with $\ones^T x= 1$, where $x_j$ is
the fraction of our wealth we place on bet $j$.
We seek to maximize the average log return, which maximizes long-term 
wealth growth if we repeatedly bet.
This leads to the (convex) optimization problem
\begin{equation*}
\begin{array}{ll}
\mbox{maximize} & \sum_{i=1}^N \pi_i \log (r_i^T x) \\
\mbox{subject to} & x \geq 0, \quad \ones^T x = 1,
\end{array}
\end{equation*}
with variable $x$.

\paragraph{Problem instance.}    
We consider $n=200$ bets. 
The probabilities of the outcomes $\pi_i$ are independent draws 
from a uniform distribution on $[0,1]$, normalized to sum to one.
The returns of the bets in each outcome are independently drawn 
from a log normal distribution, 
\ie, $\exp(z)$, $z \sim \mathcal{N}(0, 1)$, 
and then scaled so the expected return
of bet $j$ is $\bar r_j$, \ie, $\sum_{i=1}^N \pi_i (r_i)_j = \bar r_j$,
where $\bar r_j$ is drawn from a uniform distribution
on $[0.9, 1.1]$.
For our problem instance, the solution has 57
nonzero entries ranging from 0.001 to 0.04, 
and the optimal mean log return is 0.057.

\paragraph{Results.}  
The run-times for OSMM, MOSEK, SCS, 
and ECOS are shown in table~\ref{tab:runtimes:kelly}.  
When the number of Monte Carlo samples is small, \eg, $N = 1{,}000$, 
MOSEK performs the fastest and takes less than a second to attain high accuracy.
OSMM also takes less than one second.
ECOS takes about two seconds,
and SCS takes roughly six seconds.  
However, when $N$ is larger (\ie, $N = 10{,}000$, $100{,}000$, 
or $1{,}000{,}000$), OSMM is the fastest method, always taking less than a second to attain the required $10^{-6}$ optimality gap.
When $N=10{,}000$, MOSEK is still competitive, but SCS and ECOS 
are two to four orders of magnitude slower than OSMM.
When $N=100{,}000$, MOSEK and ECOS are two to three 
orders of magnitude slower, and when $N=1{,}000{,}000$, 
MOSEK and ECOS are three to five orders of magnitude slower.
SCS fails for both the two larger values of $N$.
These findings suggest that OSMM is useful when the number of samples is large, 
as it exhibits good scaling with $N$.

OSMM spends 0.0013 and 0.0024 seconds to evaluate $f$ 
and its gradient $\nabla f$, respectively, when $N=1{,}000{,}000$.  
Computing the tentative iterate $x^{k+1/2}$ and the lower bound 
$L_k$ from~\eqref{e-Lk} takes 0.033 and 0.014 seconds, respectively.  
The line search also turns out to be quite efficient here, 
as $f$ is evaluated twice during the line search on average.

\begin{table}
	\caption{Solve times in seconds on the Kelly gambling problem.  
		A dash (``---'') means the solver failed, either for numerical reasons,
		or because it did not reach the required $10^{-6}$ suboptimality
		in 24 hours.}
	\label{tab:runtimes:kelly}
	\centering
	\begin{tabular}{l llll}
		\hline\noalign{\smallskip}
			$N$ & OSMM & MOSEK & SCS & ECOS \\
		\noalign{\smallskip}\hline\noalign{\smallskip}
			1{,}000   & 0.76 & 0.58&6.4& 2.2   \\
			10{,}000 & 0.64 & 4.5 & 2{,}100 & 50 \\
			100{,}000 & 0.62 & 62 & --- & 910 \\
			1{,}000{,}000 & 0.64 & 890 & --- & 20{,}000 \\
		\noalign{\smallskip}\hline
	\end{tabular}
\end{table}

\begin{figure*}
	\centering
	\includegraphics[width=0.75\textwidth]{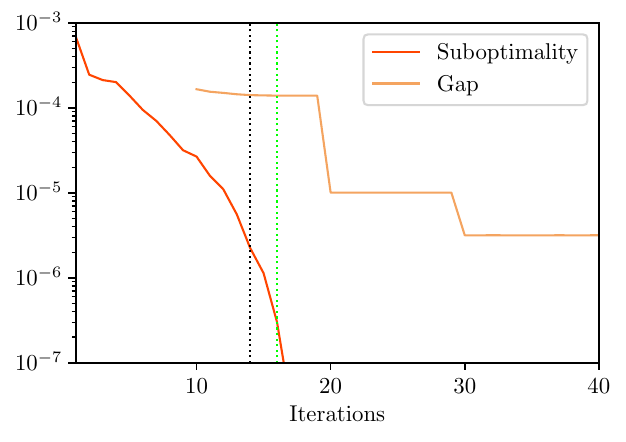}
	\includegraphics[width=0.75\textwidth]{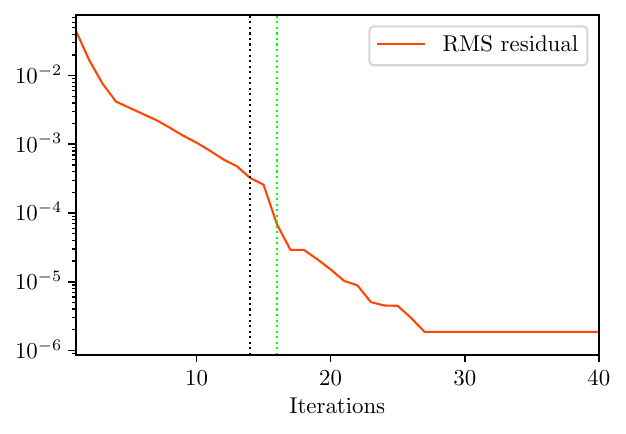}
	\caption{Suboptimality, gap (top row), 
		and RMS residual (bottom row) on the Kelly gambling problem with $N = 1{,}000{,}000$.
		The dotted green and black lines show when high accuracy and 
		practical accuracy are reached, respectively.}
	\label{fig:kelly_convergence} 
\end{figure*}

Figure~\ref{fig:kelly_convergence} shows the convergence of OSMM 
with $N = 1{,}000{,}000$.  
In the top panel, we can see that practical accuracy is reached after
14 iterations, and high accuracy is reached after 16 iterations,
as shown by the dotted black and green lines, respectively.

\subsection{CVaR portfolio optimization with derivatives}
\label{subsec:experiment-cvar}
\paragraph{Problem formulation.}
We consider making investments in $m$ stocks, and call and put
derivatives on them, with various strike prices.
Our investment will be for one period (of, say, one month).
We let $\omega \in \reals_{++}^m$ denote the (fractional) change in prices
of the $m$ underlying stocks, which we model as random with a log-normal
distribution, \ie, $\log \omega \sim \mathcal N(\mu,\Sigma)$,
where the $\log$ is elementwise.
For simplicity we will assume there is one call and one put option available
for each stock.
Let $p_\text{c} \in\reals_{++}^m$ and $s_\text{c}\in\reals_{++}^m$ 
be the call option prices (premia) and strike prices, normalized by the current 
stock price, so, for example, $(s_\text{c})_i = 1.15$ means the strike price
is $1.15$ times the current stock price.
Let $p_\text{p}\in\reals^m_{++}$ and $s_\text{p}\in\reals^m_{++}$ 
denote the corresponding quantities for the $m$ put options.

The amount we receive per dollar of investment in the underlying stocks
is $\omega$, the ratio of the current to final stock price.
For every dollar invested in the call options we receive
$(\omega-s_\text{c})_+ / p_\text{c}$, where the division is 
elementwise, and $(a)_+ = \max\{a,0\}$.
For the put options, the total we receive per dollar of investment is
$(s_\text{p}-\omega)_+ / p_\text{p}$.

We make investments in $n=3m$ different
assets, the $m$ underlying stocks and $m$ associated call and put options.
We let $x \in \reals^n$ denote the fractions of our wealth that we invest in
the assets, so $\ones^T x =1$.  We consider long and short positions, with
$x_i <0$ denoting a short position.
We consider a simple set of portfolio constraints,
$x \geq x_\text{min}$ (\eg, $x_\text{min}=-0.1$ 
limits the maximum short position for 
any asset to not exceed 10\% of the total portfolio value),
and $\|x\|_1 \leq L$, where $L$ is a leverage limit.
(Since $\ones^T x =1$, this means that the total long position cannot exceed
a multiple $(L+1)/2$ of total wealth, and the total short position
cannot exceed a multiple $(L-1)/2$ of the total wealth.)
We partition $x$ as $x=(x_\text{u},x_\text{c},x_\text{p})$, with
each subvector in $\reals^m$.
Our portfolio has total return
\[
x_\text{u}^T \omega + x_\text{c}^T ((s_\text{c}-\omega)_+ / p_\text{c} ) 
+ x_\text{p}^T ((s_\text{p}-\omega)_+ / p_\text{p}) = r(\omega)^T x,
\]
where $r(\omega) \in \reals^n$ is the total return of the $n=3m$ assets,
\ie, stocks and options.

Our problem is to choose the portfolio so as to minimize the conditional
value at risk (described in \eqref{e-cvar}) of the negative total return, \ie,
\[
\begin{array}{ll}
\mbox{minimize} &  \CVaR(-r(\omega)^T x;\eta)\\
\mbox{subject to} & x \geq x_\text{min}, \; \ones^T x = 1, \; \|x\|_1 \leq L,
\end{array}
\]
where $\eta \in (0,1)$ sets the risk aversion.

We use a sample average approximation of the expectation in CVaR to 
obtain the problem
\[
\begin{array}{ll}
\mbox{minimize} &  
\frac{1}{N} \sum_{i=1}^N
\frac{(-r(\omega_i)^Tx - \alpha)_+}{1-\eta} + \alpha \\
\mbox{subject to} & x \geq x_\text{min}, \; \ones^T x = 1, \; \|x\|_1 \leq L
\end{array}
\]
with variables $x = (x_\text{u},x_\text{c},x_\text{p}) \in \reals^n$
and $\alpha \in \reals$.
The vectors $\omega_i \in \reals_{++}^n$, $i=1,\ldots,N$, are independent 
samples of the price change $\omega$.

\paragraph{Problem instance.} 
We take the number of stocks as $m=100$, so $n=300$.
We take the minimum position limit $x_\text{min} = -0.1$
and leverage limit $L=1.6$.  We use risk aversion parameter $\eta = 0.8$,
so we are attempting to minimize the 20th percentile of the portfolio loss.
The price change covariance $\Sigma$ is generated according to 
$\Sigma = \sigma^2(I + 0.2 FF^T)$, where $\sigma = 1/\sqrt{2}$, 
and the entries of $F\in\reals^{m \times 5}$ are 
independent draws from a standard normal distribution.
We set the mean price change according to $\mu_i = 0.03 \sqrt{\Sigma_{ii}} - 0.5\Sigma_{ii}$, $i=1,\ldots,m$. 
For each call option, the strike price is set as the 80th percentile
of $\omega$, and for each put option it is the 20th percentile.
The option prices are determined by the Black-Scholes formula
with zero discount.
(These data are approximately consistent with an investment period of 
one month for U.S.~equities.)

When we solve this problem instance,
the optimal portfolio return has mean 1.3\%,
standard deviation 5\%, and loss probability 39\%;
annualized, these correspond to a return of 16\%, 
standard deviation 18\%, and loss probability 20\%.
The optimal portfolio contains as assets the underlying stocks as well as call and put options.

\begin{figure*}
	\centering
	\includegraphics[width=0.75\textwidth]{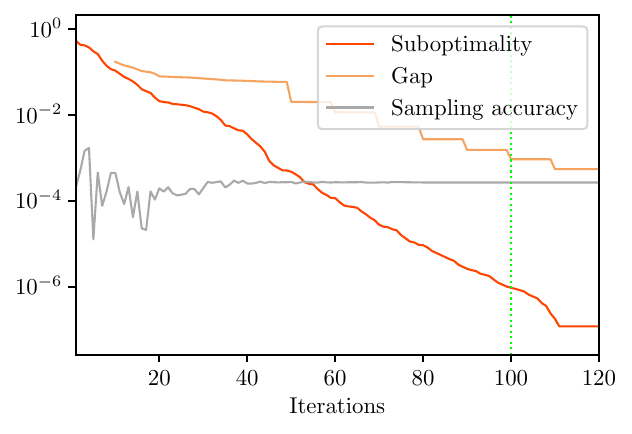}
	\includegraphics[width=0.75\textwidth]{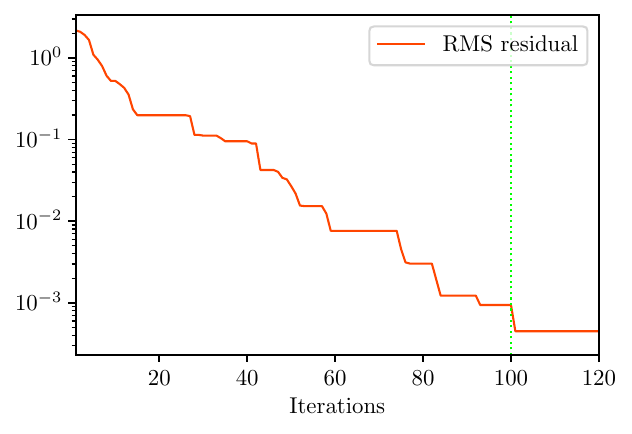}
	\caption{Suboptimality, gap, sampling accuracy (top row), and RMS residual (bottom row) 
		on the conditional value-at-risk problem with $N=1{,}000{,}000$.
		The dotted green and black lines coincide, indicating that
		high accuracy and practical accuracy are reached at the same time.}
	\label{fig:cvar_convergence} 
\end{figure*}

\paragraph{Results.}
Table~\ref{tab:runtimes:cvar} shows the run-times for
$N$ ranging from $10{,}000$ to $1{,}000{,}000$.  
We see again that for small values of $N$, it is more efficient to
solve the problem directly using a structured solver, whereas for 
large values, OSMM is far more efficient.
When $N=1{,}000{,}000$, OSMM (using PyTorch) takes 0.0021 seconds to evaluate $f$ 
and 0.0053 seconds to evaluate $\nabla f$; 
OSMM takes 0.050 seconds to compute the tentative iterate $x^{k+1/2}$, 
and 0.021 seconds to evaluate the lower bound $L_k$ (using CVXPY).
Figure~\ref{fig:cvar_convergence} shows the convergence of OSMM 
with $N=1{,}000{,}000$.  
Practical accuracy are high accuracy are reached at the same time
after 100 iterations.

\begin{table}
	\caption{Solve times in seconds on the conditional value-at-risk problem.  
		A dash (``---'') means the solver failed, either for numerical reasons,
		or because it did not reach the required $10^{-6}$ suboptimality
		in 24 hours.}
	\label{tab:runtimes:cvar}
	\centering
	\begin{tabular}{l llll}
		\hline\noalign{\smallskip}
			$N$ & OSMM & MOSEK & SCS & ECOS \\
		\noalign{\smallskip}\hline\noalign{\smallskip}
			10{,}000 & 11 & 6.0 & 250 & 18 \\
			100{,}000 & 6.5 & 64 & 2{,}900 & 310\\
			1{,}000{,}000 & 6.3 & 1{,}900 & 30{,}000 & 5{,}200 \\ 
		\noalign{\smallskip}\hline
	\end{tabular}
\end{table}

\subsection{Exponential series density estimation}\label{subsec:experiment:ESE}
\paragraph{Problem formulation.} We consider an instance of the 
generic exponential family density fitting problem 
described in \S\ref{s-exp-fam-fit}.
We consider data $z_1, \ldots, z_m \in \reals^2$, 
and wish to fit a density $p$ supported on $\mathcal S = [-1,1]^2$.
We let the sufficient statistics $\phi_i, \; i=1,\ldots,n$, be 
the Legendre polynomials up to degree four, so $n=14$.
(This is known as an exponential series 
density estimator~\cite{marsh2007goodness,wu2010exponential,gao2015penalized}.)
We solve the density estimation problem~\eqref{exp_density_prob_0}.

\paragraph{Problem instance.}   We take $m$ data points
sampled from a mixture of three Gaussian densities,
restricted to $\mathcal S$, with means $(1/3,1/3)$, $(1/3,-1/3)$, and
$(-1/3,-1/3)$, weights $0.4$, $0.3$, and $0.3$, 
and common covariance $(1/36) I$.
(So the data do not come from the family of density we use to fit.)
We form a Riemann sum using points in $\mathcal S$ lying on a grid with side lengths $\sqrt N$.
Recall that $N$, the number of samples, here refers to our approximate evaluation
of the integral that arises in the log-partition function, and not 
the number of data samples, which is fixed at $m=2{,}000$.

\paragraph{Results.}  
Table~\ref{tab:runtimes:density} shows the run-times for the various methods.
We see the usual pattern where directly solving the problem can be efficient for small $N$, but OSMM is much faster for large $N$.
When $N=1{,}000{,}000$,
it takes OSMM 0.001 and 0.0014 seconds to evaluate $f$ and $\nabla f$, 
respectively,
and 0.015 and 0.0097 seconds to compute $x^{k+1/2}$ and $L_k$, respectively.

\begin{figure*}
	\centering
	\includegraphics[width=0.75\textwidth]{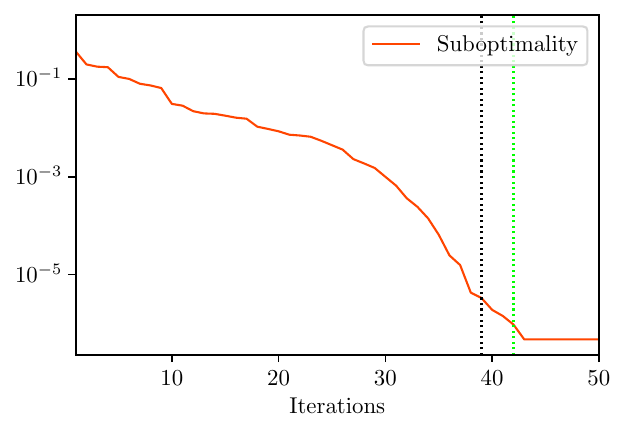}
	\includegraphics[width=0.75\textwidth]{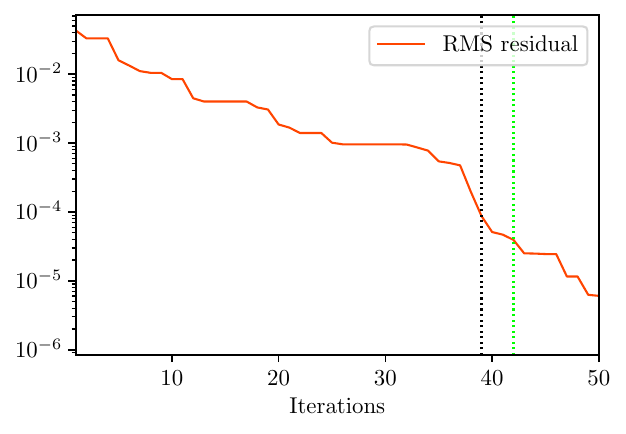}
	\caption{Suboptimality (top row) and RMS residual (bottom row) 
	on the exponential family density fitting problem with $N=1{,}000{,}000$.
	The dotted green and black lines show when high accuracy and 
	practical accuracy are reached, respectively.}
	\label{fig:exp_density_Legendre_polynomials} 
\end{figure*}

Figure~\ref{fig:exp_density_Legendre_polynomials} shows the convergence 
of OSMM for $N=1{,}000{,}000$.
Practical accuracy is reached after 39 iterations with suboptimality $10^{-5}$, 
and after 42 iterations OSMM reaches high accuracy.
In this instance, our lower bound $L_k = -\infty$, so neither it nor the gap are plotted in the figure.

\begin{table}
	% The table caption is above the table.
	\caption{Solve times in seconds on the density estimation problem.  A dash (``---'') means the solver failed, either for numerical reasons,
		or because it did not reach the required $10^{-6}$ suboptimality
		in 24 hours.}
	\label{tab:runtimes:density}
	\centering
		\begin{tabular}{l llll}
			\hline\noalign{\smallskip}
			$N$ & OSMM & MOSEK & SCS & ECOS \\
			\noalign{\smallskip}\hline\noalign{\smallskip}
			10{,}000 & 0.72 & 1.2 & 1{,}100 & 0.92 \\
			100{,}000 & 1.2 & 11 & --- & 14  \\ 
			1{,}000{,}000 &  0.84 &  120  & --- & 190 \\
			\noalign{\smallskip}\hline
		\end{tabular}
\end{table}

\subsection{Vector news vendor with entropic value-at-risk}
We consider a variant of the classic news vendor problem that 
involves entropic value-at-risk programming~\cite{ahmed2007coherent}.

\paragraph{Problem formulation.} 
We choose quantities $q \in \reals_+^n$ of $n$
products to produce, at total cost $\phi(q)$, where $\phi : \reals^n_+ \to
\reals_+$.  We have constraints on the quantities we can produce, 
$q \leq q_\mathrm{max}$, and on the total production cost, $\phi(q) \leq 
\phi_\mathrm{max}$.
We sell the amount $\min(q,d)$, where $d\in \reals_+^n$ is the demand,
and the min is taken elementwise, at market prices $p \in \reals_+^n$,
so the total revenue is $p^T \min(q,d)$, and the profit is
\[
P = p^T \min\{q,d\} - \phi(q).
\]
We assume that $\omega = (d,p) \in \reals_+^{2n}$ is a random variable
with known distribution, so $P(q,\omega)$ 
is a random variable that depends on $q$.

We choose the quantities $q$ to minimize the EVaR of the negative
profit,
\[
\begin{array}{ll}
\mbox{minimize} & \EVaR (-P(q,\omega); \eta) \\ 
\mbox{subject to} &  \phi(q) \leq \phi_{\max}, \; q \geq 0, \; q \leq q_{\max},
\; \end{array}
\]
where $\eta$ is a specified quantile.

As described in \S\ref{s-EVaR}, 
we approximate this with Monte Carlo samples $(d_1,p_1), \ldots, 
(d_N,p_N)$ to obtain the problem
\[
\begin{array}{ll}
\mbox{minimize} & \alpha \log\left(\frac{1}{(1-\eta)N} \sum_{i=1}^N \exp
\left(\frac{-p_i^T\min(q, d_i) + \phi(q)}{\alpha} \right) \right)\\
\mbox{subject to} &  \phi(q) \leq \phi_{\max}, \; q \geq 0, \; q \leq q_{\max},
\; \alpha \geq 0, \end{array}
\]
with variables $q$ and $\alpha\in\reals_+$.  We
denote the realizations of the demand $d$ and prices $p$ on the $i$th Monte
Carlo simulation by $d_i, p_i \in \reals_+^n, \; i=1,\ldots,N$, respectively.

\paragraph{Problem instance.}  
We take $n=500$ and risk aversion parameter $\eta=0.9$.
We assume the demand and market prices follow 
a joint log-normal distribution, \ie, 
$(d,p) = \exp z$, $z \sim \mathcal N(\mu, \Sigma)$, 
and the exponential is elementwise.  
We draw the entries of $\mu\in \reals^{2n}$ independently from 
a uniform distribution on $[-0.2, 0]$,
and set $\Sigma = 0.1FF^T$,
where the entries of $F\in \reals^{2n \times 5}$ are independently 
drawn from a standard normal distribution. 
We use a production cost which is separable and piecewise affine
with one kink point for each entry of $q$,
\[
\phi(q) = a^T q + 0.5 a^T (q - b)_+,
\]
where the elements in $a$ and $b$ are both drawn uniformly at random from $[0.2, 0.9]$ and $[0.01, 0.03]$, respectively.
The maximum production quantities $q_{\max}$ is set as 
$5b$, and the maximum cost is $\phi_\textrm{max} = 1$.
With these parameter values, the optimal
profit has mean 3.2 and standard deviation 0.98.

\paragraph{Results.}  
The run-times for the various methods are in 
table~\ref{tab:runtimes:vecnewsven}.
In this instance, OSMM is the fastest for all values of $N$ ranging from 
$1{,}000$ to $1{,}000{,}000$, and the other solvers fail for nearly all values of $N$.
When $N=1{,}000{,}000$, OSMM takes 0.16 and 0.27 seconds to evaluate 
$f$ and $\nabla f$, respectively; it takes 0.076 seconds to compute 
the tentative iterate, and 0.036 seconds to compute the lower bound. 

\begin{figure*}
	\centering
	\includegraphics[width=0.75\textwidth]{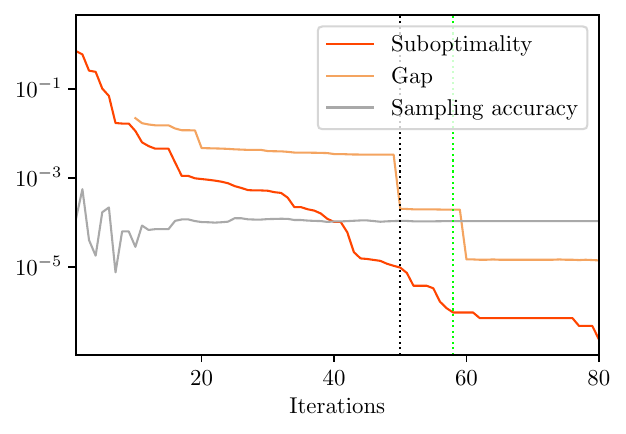}
	\includegraphics[width=0.75\textwidth]{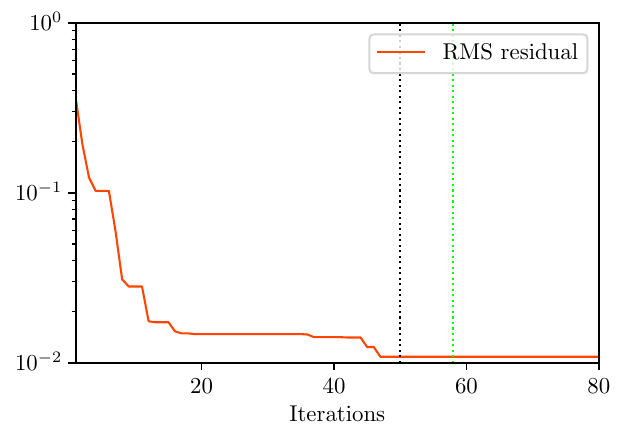}
	\caption{Suboptimality, gap, sampling accuracy (top row), and RMS residual (bottom row) 
		on the vector news vendor problem with $N=1{,}000{,}000$.
		The dotted green and black lines show when high accuracy and 
		practical accuracy are reached, respectively.}
	\label{fig:vecnewsven_convergence} 
\end{figure*}

Figure~\ref{fig:vecnewsven_convergence} shows 
the convergence of OSMM with $N=1{,}000{,}000$.  
We can see that practical accuracy is reached after 50 
iterations with suboptimality on the order of $10^{-5}$, 
while it takes 58 iterations to reach high accuracy.

\begin{table}
	% The table caption is above the table.
	\caption{Solve times in seconds on the vector news vendor problem.  
		A dash (``---'') means the solver failed, either for numerical reasons,
		or because it did not reach the required $10^{-6}$ suboptimality
		in 24 hours.}
	\label{tab:runtimes:vecnewsven}
	\centering
		\begin{tabular}{l llll}
			\hline\noalign{\smallskip}
			$N$ & OSMM & MOSEK & SCS & ECOS \\
			\noalign{\smallskip}\hline\noalign{\smallskip}
			1{,}000 & 6.7 & 120 &  --- & ---\\
			10{,}000 & 11 & 7{,}400 & --- & --- \\
			100{,}000 & 10 & --- & --- & ---  \\ 
			1{,}000{,}000 & 53  &  ---  & --- & --- \\
			\noalign{\smallskip}\hline
		\end{tabular}
\end{table}

\appendix  %This command ends the counting of sections.
\section{Appendix}

\subsection{Details of forming $G_k$} \label{appendix:s-curvature}
We form $G_k$ in \eqref{PSD_low_rank}
by adopting a quasi-Newton update given in \cite{fletcher2005new}.
The main idea is to divide $G_k$ into two matrices
$G_k^{(1)}$ and $G_k^{(2)}$, which are the first $r_1$
and the last $r_2=r-r_1$ columns in $G_k$, respectively,
and update them separately into
$G_{k+1}^{(1)}$ and $G_{k+1}^{(2)}$. Then $G_{k+1}$ is updated by
\[
G_{k+1} = \left[G^{(1)}_{k+1} \quad G^{(2)}_{k+1}\right].
\]

The detail is as follows.
Let $s_k = x^{k+1} - x^k$ and $y_k = \nabla f(x^{k+1}) - \nabla f(x^k)$. 
From the convexity of $f$, $s_k^T y_k\geq 0$.
Suppose that $s_k^T y_k$ is not too small such that
\BEQ\label{s_ky_k_bound}
s_k^T y_k > \max(\varepsilon_\mathrm{abs},
\varepsilon_\mathrm{rel}\|s_k\|_2\|y_k\|_2),
\EEQ
where constants $\varepsilon_\mathrm{abs}, \varepsilon_\mathrm{rel}>0$ 
are given.
Then $r_1$ is chosen as the largest integer in $[0, r]$ such that 
$G_k^{(1)}$ satisfies
\[
s_k^Ty_k - \left\|(G_k^{(1)})^Ts_k\right\|_2^2 > 
\varepsilon_\mathrm{rel} \|s_k\|_2 \left\|y_k - G_k^{(1)}(G_k^{(1)})^Ts_k \right\|_2.
\]
According to \eqref{s_ky_k_bound} the above holds
at least for $r_1=0$, in which case
$G^{(1)}_k$ degenerates to $0$.

Once $r_1$ is obtained, we define
$w_k^{(1)} = (G^{(1)}_k)^T s_k$ and
$w_k^{(2)} = (G^{(2)}_k)^T s_k$.
Then $R_k^{(1)}\in\reals^{(r_1+1)\times(r_1+1)}$ is the upper triangular factor in the
following R-Q decomposition
\begin{align*}%\label{rq_fac}
R_k^{(1)} Q_k^{(1)} = 
\left[
\begin{array}{ll}
\frac{1}{\sqrt{s_k^T y_k - \|w_k^{(1)}\|_2^2}} & 0_{r_1}^T \\
\frac{-1}{\sqrt{s_k^T y_k - \|w_k^{(1)}\|_2^2}}w_k^{(1)} & I_{r_1}
\end{array}
\right] \in\reals^{(r_1+1)\times(r_1+1)},
\end{align*}
and $Q_k^{(2)}\in\reals^{r_2 \times (r_2 - 1)}$ is a set of 
orthonormal basis orthogonal to $w_k^{(2)}$.
The update is
\begin{align*}%\label{update_U}
G^{(1)}_{k+1} = \left[ y_k \quad G^{(1)}_k \right]R_k^{(1)}, \quad
G^{(2)}_{k+1} = G^{(2)}_{k} Q_k^{(2)}.
\end{align*} 

There are some corner cases.
If \eqref{s_ky_k_bound} holds and $r_1 = 0$, 
then $G_{k+1}^{(1)} = y_k / \sqrt{s_k^T y_k}$.
If $r_1 = r-1$ or $r_1 = r$, then $G_{k+1}^{(2)}$ vanishes,
and $G_{k+1}$ takes the first $r$ columns of $G_{k+1}^{(1)}$.

In cases where \eqref{s_ky_k_bound} does not hold,
if $\|w_k^{(2)}\|_2 > \varepsilon_\mathrm{abs}$,
then we can still define $G_{k+1}^{(2)}$ in the same way,
and $G_{k+1} = \left[G^{(2)}_{k+1} \quad 0_n\right].$
Otherwise, $G_{k+1} = G_k$.

It can be easily checked that by the $G_k$ defined as above,
the trace of $H_k$ is uniformly upper bounded in $k$.

The default values for the parameters are 
$\varepsilon_\mathrm{abs}=10^{-8}$ and $\varepsilon_\mathrm{rel}=10^{-3}$.

\subsection{Details of computing an optimal subgradient of $g$} \label{appendix-q_k}

Here we show how to compute a subgradient $q^{k+1} \in \partial g(x^{k+1/2})$ satisfying the optimality conditions in \eqref{e-xk+half-opt}, which implies \eqref{e-qk}.  This, in turn, is important because it allows us to compute the stopping criteria described in \S\ref{s-stopping-crit}.

Since we know the third term on the right-hand side of \eqref{e-xk+half-opt}, it suffices to find an optimal subgradient $l_k^{k+1/2} \in \partial l_k(x^{k+1/2})$, which is easier.  We start by rewriting the defining problem for $x^{k+1/2}$ in a more useful form.  Putting \eqref{l_k_PWA} and \eqref{e-xk+half2} together, we see that we can alternatively express $x^{k+1/2}$ as the solution to the convex problem
\BEQ\label{e-q-prob}
\begin{array}{ll}
	\mbox{minimize} & z + g(x) + \frac{1}{2}(x-x^k)^T (H_k + \lambda_k I) (x-x^k) \\
	\mbox{subject to} & z \geq f(x^i) + \nabla f(x^i)^T (x - x^i), 
	\quad \textrm{for } i=\max\{0, k-M+1\},\ldots,k,
\end{array}
\EEQ
with variables $x$ and $z \in \reals$.

The KKT conditions for problem \eqref{e-q-prob}, which are necessary and sufficient for the points $(x^{k+1/2}, z^{k+1/2})$ and $\gamma_i, \; i=\max\{0, k-M+1\},\ldots,k$, to be primal and dual optimal, are
\begin{align}
z^{k+1/2} \geq f(x^i) + \nabla f(x^i)^T (x^{k+1/2} - x^i), \quad & i=\max\{0, k-M+1\},\ldots,k; \label{KKT_3} \\ %\quad \textrm{(primal feasibility and stationarity w.r.t.~$z$)} 
\sum_{i=\max\{0, k-M+1\}}^k \gamma_i = 1, \quad \gamma_i \geq 0, \quad & i=\max\{0, k-M+1\},\ldots,k; \label{KKT_1} \\ %\quad \textrm{(dual feasibility)} 
z^{k+1/2} > f(x^i) + \nabla f(x^i)^T (x^{k+1/2} - x^i) \implies \gamma_i = 0, \quad & i\in\max\{0, k-M+1\},\ldots,k; \label{KKT_4} \\ %\quad \textrm{(complementary slackness)}
\partial g(x^{k+1/2}) + (H_k + \lambda_k I) v^k 
+ \sum_{i=\max\{0, k-M+1\}}^k 
\gamma_i \nabla f(x^i) \ni & 0. \label{KKT_2} % \quad \textrm{(stationarity w.r.t.~$x$).} 
\end{align}
Here we used the definition of $v^k$ to simplify the stationarity condition \eqref{KKT_2}.

Now we claim that
\[
l_k^{k+1/2} = \sum_{i=\max\{0, k-M+1\}}^k \gamma_i \nabla f(x^i) \in \partial l_k(x^{k+1/2}).
\]
To see this, note that \eqref{KKT_1} says the $\gamma_i$ are nonnegative and sum to one, while \eqref{KKT_3} and \eqref{KKT_4} together imply $\gamma_i$ is positive as long as $\nabla f(x^i)$ is active; this satisfies \eqref{e-subdiff-lk}, which says the subdifferential $\partial l_k(x^{k+1/2})$ is the convex hull of the active gradients.  Therefore, re-arranging \eqref{KKT_2} gives
\[
q^{k+1} = - l_k^{k+1/2} - (H_k + \lambda_k I) v^k \in \partial g(x^{k+1/2}),
\]
which shows \eqref{e-qk}.

\subsection{Proof that undamped steps occur infinitely often}\label{s-tk-one-io}
Assume that $\nabla f$ is Lipschitz continuous with constant $L$.  Also, assume $\mu_{\max}\tau_{\min} > 2L / (1-\alpha)$.
We show that for any number of iterations $k_0$, there is some $k \geq k_0$ such that $t_k = 1$.  This means that there exists a subsequence $(k_\ell)_{\ell=1}^\infty$ such that $t_{k_\ell} = 1$.

To show the result, we first claim that the line search condition \eqref{armijo_step} is satisfied as soon as
\BEQ\label{t_k_upper_bound}
t_k \leq \frac{1-\alpha}{2L} \lambda_k.
\EEQ
We prove the claim later.  Taking the claim as a given for now, the main result follows by deriving a contradiction.  
To get a contradiction, suppose there exists some number of iterations $k_0$ such that $t_k < 1$ for each $k \geq k_0$.  
Then, because $t_k$ is the largest number satisfying 
	$t_k=\beta^j$ ($j \in \mbox{\bf{N}}_0$) and the line search condition,
	we get that the line search condition does not hold with $t_k=1$,
	and thus \eqref{t_k_upper_bound} does not hold with $t_k=1$, 
	meaning that 
$\lambda_k < 2 L / (1-\alpha)$, for each $k \geq k_0$.  But from \eqref{e-muk}, we have $\mu_k = \min\left\{\gamma_{\mathrm{inc}}^{k-k_0} \mu_{k_0},\mu_{\max}\right\}$, since we assumed $t_k < 1$ for every $k \geq k_0$.  Additionally, from \eqref{e-lamk}, we get $\lambda_k \geq \mu_k \tau_{\min}$.  So, we now have two cases.  Either we have $\lambda_k \geq \mu_{\max} \tau_{\min} > 2L / (1-\alpha)$, 
which is a contradiction with $\lambda_k < 2 L / (1-\alpha)$.  
Or we have $\lambda_k \geq \gamma_{\mathrm{inc}}^{k-k_0} \mu_{k_0}\tau_{\min}$, in which case $\gamma_{\mathrm{inc}}^{k-k_0} \mu_{k_0}\tau_{\min}$ grows exponentially in $k$; this means we must have $\lambda_k \geq 2 L / (1-\alpha)$, for $k$ sufficiently large, which is again a contradiction.  This finishes the proof of the main result.

Now we prove the claim.  Observe that
\BEQ\label{e-t_k-upper-bound-2}
L t_k^2 \|v^k\|^2 \leq \frac{1-\alpha}{2} \lambda_k t_k \|v^k\|^2
\leq \frac{1-\alpha}{2} t_k (v^k)^T (H_k +\lambda_k I) v^k,
\EEQ
where we used \eqref{t_k_upper_bound} to get the first inequality.  We now use the following two facts: $\nabla f$ being $L$-Lipschitz continuous implies that (i) $\phi_k$ is convex in $t$, and (ii) $\phi_k^\prime$ is $L \| v^k \|^2$-Lipschitz in $t$. %, so we get 
By the first fact (convexity), we have
\[
\phi_k(t_k) \leq \phi_k(0) + \phi_k^\prime(t_k) t_k.
\]
Adding and subtracting $\phi_k^\prime(0)$ gives
\[
\phi_k(t_k) \leq \phi_k(0) + \phi_k^\prime(0) t_k + (\phi_k^\prime(t_k) - \phi_k^\prime(0)) t_k.
\]
The second fact (Lipschitz continuity) yields a bound on the third term on the right-hand side,
\BEQ\label{e-phi-tk-upper-bound}
\phi_k(t_k) \leq \phi_k(0) + \phi_k^\prime(0) t_k + L t_k^2 \|v^k\|^2.
\EEQ
Finally, using \eqref{phi_prime_bound} and \eqref{e-t_k-upper-bound-2} to bound $\phi_k^\prime(0)$ and $L t_k^2 \| v^k \|_2^2$ in \eqref{e-phi-tk-upper-bound} yields \eqref{armijo_step}, which implies the line search condition \eqref{armijo_step} is satisfied, as claimed.

We note that the claim above also implies the following lower bound, which is used in \S\ref{s-convergence}.   
Observe that if $t_k$ indeed satisfies the line search condition, then we can express $t_k = \beta^j$, where $j$ is the smallest integer such that \eqref{t_k_upper_bound} holds (see \S\ref{s-line-search}).  Now consider two cases.  If $1 \leq (1-\alpha) \lambda_k / (2L)$, then we can simply take $j = 0$, so that $t_k = 1$.  On the other hand, if $1 > (1-\alpha) \lambda_k / (2L)$, then a short calculation shows that $j=\lceil \log_{\beta} ((1-\alpha) \lambda_k / (2L)) \rceil$, and so $j < 1 + \log_{\beta} ((1-\alpha) \lambda_k / (2L))$, giving $t_k = \beta^j > \beta (1-\alpha) \lambda_k / (2L)$.  Therefore, to sum up, we have \eqref{t_k_upper_bound} implies that $t_k$ satisfies the inequalities
\BEQ\label{lower-bound-tk}
t_k > \beta \min\left\{1, \frac{1-\alpha}{2L} \lambda_k \right\}
\geq  \beta \min\left\{1,\frac{1-\alpha}{2L} \mu_{\min}\tau_{\min}\right\},
\EEQ
where we also used the fact that $\lambda_k \geq \mu_{\min}\tau_{\min}$.

\subsection{Proof that the limited memory piecewise affine minorant is accurate enough}\label{s-minorant-accuracy}

Assume that $\nabla f$ is Lipschitz continuous with constant $L$.  For the rest of the proof, fix any $k$ such that $t_k = 1$.  (From \S\ref{s-tk-one-io}, it is always possible to find at least one such $k$.)  We will show that for any such $k$, the limited memory piecewise affine minorant \eqref{l_k_PWA} satisfies the bound \eqref{e-minorant-accuracy}; because our choice of $k$ was arbitrary, the required result will then follow.

For any $l_k^{k+1} \in \partial l_k(x^{k+1})$, note that adding and subtracting $\nabla f(x^k)$ in the left-hand side of \eqref{e-minorant-accuracy} easily gives
\BEQ\label{e-minorant-accuracy-proof}
\| \nabla f(x^{k+1}) - l_k^{k+1} \|_2 \leq \| \nabla f(x^{k+1}) - \nabla f(x^k) \|_2 + \| \nabla f(x^k) - l_k^{k+1} \|_2.
\EEQ
The Lipschitz continuity of $\nabla f$, in turn, immediately gives a bound for the first term on the right-hand side of \eqref{e-minorant-accuracy-proof}, \ie, we get
\begin{align}
\| \nabla f(x^{k+1}) - l_k^{k+1} \|_2 & \leq L \| x^{k+1} - x^k \|_2 + \| \nabla f(x^k) - l_k^{k+1} \|_2 \notag \\
& \lesssim \| v^k \|_2 + \| \nabla f(x^k) - l_k^{k+1} \|_2, \label{e-minorant-accuracy-4-proof}
\end{align}
using the definition of $v^k$.

Therefore, we focus on the second term on the right-hand side of \eqref{e-minorant-accuracy-proof}.  For this term, \eqref{e-subdiff-lk} tells us that for any $l_k^{k+1} \in \partial l_k(x^{k+1})$,
\BEQ\label{e-minorant-accuracy-2-proof}
\| \nabla f(x^k) - l_k^{k+1} \|_2 \leq \max_{j=\max\{0,k-M+1\},\ldots,k} \| \nabla f(x^k) - \nabla f(x^j) \|_2,
\EEQ
because the maximum of a convex function over a convex polytope is attained at one of its vertices.  The Lipschitz continuity of $\nabla f$ again shows that, for any $j \in \{ \max\{0,k-M+1\}, \ldots, k \}$,
\begin{align*}
\| \nabla f(x^k) - \nabla f(x^j) \|_2 & \leq L \| x^k - x^j \|_2 \\
& = L \| (x^k - x^{k-1}) + (x^{k-1} - x^{k-2}) + \cdots + (x^{j+2} - x^{j+1}) + (x^{j+1} - x^j) \|_2 \\
& \leq L \sum_{\ell=j}^{k-1} \| v^\ell \|_2 \\
& \lesssim \max_{\ell=j,\ldots,k-1} \| v^\ell \|_2.
\end{align*}
To get the third line, we used the fact that the sum on the second line telescopes, and applied the triangle inequality.  To get the fourth line, we used the fact that the average is less than the max.  Putting this last inequality together with \eqref{e-minorant-accuracy-2-proof}, we see that
\BEQ\label{e-minorant-accuracy-3-proof}
\| \nabla f(x^k) - l_k^{k+1} \|_2 \lesssim \max_{j=\max\{0,k-M+1\},\ldots,k} \| v^j \|_2,
\EEQ
for any $l_k^{k+1} \in \partial l_k(x^{k+1})$.  Combining \eqref{e-minorant-accuracy-3-proof} and \eqref{e-minorant-accuracy-4-proof} gives \eqref{e-minorant-accuracy}.  This completes the proof of the result.

%\begin{acknowledgements}
%\end{acknowledgements}


\bibliographystyle{spmpsci}
\bibliography{refs}

\begin{thebibliography}{10}
\providecommand{\url}[1]{{#1}}
\providecommand{\urlprefix}{URL }
\expandafter\ifx\csname urlstyle\endcsname\relax
  \providecommand{\doi}[1]{DOI~\discretionary{}{}{}#1}\else
  \providecommand{\doi}{DOI~\discretionary{}{}{}\begingroup
  \urlstyle{rm}\Url}\fi

\bibitem{tensorflow2015-whitepaper}
Abadi, M., Agarwal, A., Barham, P., Brevdo, E., Chen, Z., Citro, C., Corrado,
  G.S., Davis, A., Dean, J., Devin, M., Ghemawat, S., Goodfellow, I., Harp, A.,
  Irving, G., Isard, M., Jia, Y., Jozefowicz, R., Kaiser, L., Kudlur, M.,
  Levenberg, J., Man\'{e}, D., Monga, R., Moore, S., Murray, D., Olah, C.,
  Schuster, M., Shlens, J., Steiner, B., Sutskever, I., Talwar, K., Tucker, P.,
  Vanhoucke, V., Vasudevan, V., Vi\'{e}gas, F., Vinyals, O., Warden, P.,
  Wattenberg, M., Wicke, M., Yu, Y., Zheng, X.: {TensorFlow}: Large-scale
  machine learning on heterogeneous systems (2015).
\newblock \urlprefix\url{https://www.tensorflow.org/}.
\newblock Software available from tensorflow.org

\bibitem{van2016strongly}
van Ackooij, W., Bello~Cruz, J.Y., de~Oliveira, W.: A strongly convergent
  proximal bundle method for convex minimization in {H}ilbert spaces.
\newblock Optimization \textbf{65}(1), 145--167 (2016)

\bibitem{adamson1969slang}
Adamson, D.S., Winant, C.W.: A {SLANG} simulation of an initially strong shock
  wave downstream of an infinite area change.
\newblock In: Proceedings of the Conference on Applications of
  Continuous-System Simulation Languages, pp. 231--240 (1969)

\bibitem{agrawal2018rewriting}
Agrawal, A., Verschueren, R., Diamond, S., Boyd, S.: A rewriting system for
  convex optimization problems.
\newblock Journal of Control and Decision \textbf{5}(1), 42--60 (2018)

\bibitem{ahmadi2011information}
Ahmadi-Javid, A.: An information-theoretic approach to constructing coherent
  risk measures.
\newblock In: 2011 IEEE International Symposium on Information Theory
  Proceedings, pp. 2125--2127. IEEE (2011)

\bibitem{ahmed2007coherent}
Ahmed, S., {\c{C}}akmak, U., Shapiro, A.: Coherent risk measures in inventory
  problems.
\newblock European Journal of Operational Research \textbf{182}(1), 226--238
  (2007)

\bibitem{mosek}
ApS, M.: MOSEK Optimizer API for Python 9.2.40 (2019).
\newblock \urlprefix\url{https://docs.mosek.com/9.2/pythonapi/index.html}

\bibitem{armijo1966minimization}
Armijo, L.: Minimization of functions having {L}ipschitz continuous first
  partial derivatives.
\newblock Pacific Journal of mathematics \textbf{16}(1), 1--3 (1966)

\bibitem{asi2019modeling}
Asi, H., Duchi, J.: Modeling simple structures and geometry for better
  stochastic optimization algorithms.
\newblock In: The 22nd International Conference on Artificial Intelligence and
  Statistics, pp. 2425--2434 (2019)

\bibitem{Bagirov2014}
Bagirov, A., Karmitsa, N., M{\"a}kel{\"a}, M.M.: Bundle Methods. Introduction
  to Nonsmooth Optimization: Theory, Practice and Software, pp. 305--310.
\newblock Springer International Publishing, Cham (2014).
\newblock \doi{10.1007/978-3-319-08114-4_12}.
\newblock \urlprefix\url{https://doi.org/10.1007/978-3-319-08114-4_12}

\bibitem{baydin2018automatic}
Baydin, A.G., Pearlmutter, B.A., Radul, A.A., Siskind, J.M.: Automatic
  differentiation in machine learning: a survey.
\newblock Journal of machine learning research \textbf{18} (2018)

\bibitem{becker2019quasi}
Becker, S., Fadili, J., Ochs, P.: On quasi-{N}ewton forward-backward splitting:
  Proximal calculus and convergence.
\newblock SIAM Journal on Optimization \textbf{29}(4), 2445--2481 (2019)

\bibitem{becker2012quasi}
Becker, S., Fadili, M.J.: A quasi-{N}ewton proximal splitting method.
\newblock arXiv preprint arXiv:1206.1156  (2012)

\bibitem{BoydVand04}
Boyd, S., Vandenberghe, L.: Convex Optimization.
\newblock Cambridge University Press (2004)

\bibitem{broyden1965class}
Broyden, C.G.: A class of methods for solving nonlinear simultaneous equations.
\newblock Mathematics of computation \textbf{19}(92), 577--593 (1965)

\bibitem{busseti2016risk}
Busseti, E., Ryu, E.K., Boyd, S.: Risk-constrained {K}elly gambling.
\newblock The Journal of Investing \textbf{25}(3), 118--134 (2016)

\bibitem{byrd1994representations}
Byrd, R.H., Nocedal, J., Schnabel, R.B.: Representations of quasi-{N}ewton
  matrices and their use in limited memory methods.
\newblock Mathematical Programming \textbf{63}(1), 129--156 (1994)

\bibitem{danielsson2013fat}
Dan{\'\i}elsson, J., Jorgensen, B.N., Samorodnitsky, G., Sarma, M., de~Vries,
  C.G.: Fat tails, var and subadditivity.
\newblock Journal of econometrics \textbf{172}(2), 283--291 (2013)

\bibitem{danielsson2008optimal}
Dan{\'\i}elsson, J., Jorgensen, B.N., de~Vries, C.G., Yang, X.: Optimal
  portfolio allocation under the probabilistic var constraint and incentives
  for financial innovation.
\newblock Annals of finance \textbf{4}(3), 345--367 (2008)

\bibitem{davidon1959variable}
Davidon, W.C.: Variable metric method for minimization.
\newblock Tech. rep., Argonne National Lab., Lemont, Ill. (1959)

\bibitem{de2016doubly}
De~Oliveira, W., Solodov, M.: A doubly stabilized bundle method for nonsmooth
  convex optimization.
\newblock Mathematical programming \textbf{156}(1-2), 125--159 (2016)

\bibitem{dennis1977quasi}
Dennis Jr, J.E., Mor{\'e}, J.J.: Quasi-{N}ewton methods, motivation and theory.
\newblock SIAM review \textbf{19}(1), 46--89 (1977)

\bibitem{cvxpy_paper}
Diamond, S., Boyd, S.: {CVXPY}: A {P}ython-embedded modeling language for
  convex optimization.
\newblock Journal of Machine Learning Research  (2016).
\newblock \urlprefix\url{http://stanford.edu/~boyd/papers/pdf/cvxpy_paper.pdf}

\bibitem{Domahidi2013ecos}
Domahidi, A., Chu, E., Boyd, S.: {ECOS}: {A}n {SOCP} solver for embedded
  systems.
\newblock In: European Control Conference (ECC), pp. 3071--3076 (2013)

\bibitem{erdogdu2015convergence}
Erdogdu, M.A., Montanari, A.: Convergence rates of sub-sampled {N}ewton
  methods.
\newblock Advances in Neural Information Processing Systems \textbf{28},
  3052--3060 (2015)

\bibitem{fletcher2005new}
Fletcher, R.: A new low rank quasi-{N}ewton update scheme for nonlinear
  programming.
\newblock In: IFIP Conference on System Modeling and Optimization, pp.
  275--293. Springer (2005)

\bibitem{fletcher1963rapidly}
Fletcher, R., Powell, M.: A rapidly convergent descent method for minimization.
\newblock The computer journal \textbf{6}(2), 163--168 (1963)

\bibitem{frangioni2002generalized}
Frangioni, A.: Generalized bundle methods.
\newblock SIAM Journal on Optimization \textbf{13}(1), 117--156 (2002)

\bibitem{fukushima1984descent}
Fukushima, M.: A descent algorithm for nonsmooth convex optimization.
\newblock Mathematical Programming \textbf{30}(2), 163--175 (1984)

\bibitem{gao2015penalized}
Gao, Y., Zhang, Y.Y., Wu, X.: Penalized exponential series estimation of copula
  densities with an application to intergenerational dependence of body mass
  index.
\newblock Empirical Economics \textbf{48}(1), 61--81 (2015)

\bibitem{ghanbari2018proximal}
Ghanbari, H., Scheinberg, K.: Proximal quasi-{N}ewton methods for regularized
  convex optimization with linear and accelerated sublinear convergence rates.
\newblock Computational Optimization and Applications \textbf{69}(3), 597--627
  (2018)

\bibitem{gower2018tracking}
Gower, R., Le~Roux, N., Bach, F.: Tracking the gradients using the hessian: A
  new look at variance reducing stochastic methods.
\newblock In: International Conference on Artificial Intelligence and
  Statistics, pp. 707--715. PMLR (2018)

\bibitem{grant2014cvx}
Grant, M., Boyd, S.: {CVX}: Matlab software for disciplined convex programming,
  version 2.1 (2014)

\bibitem{GBY:06}
Grant, M., Boyd, S., Ye, Y.: Disciplined convex programming.
\newblock In: L.~Liberti, N.~Maculan (eds.) Global Optimization: From Theory to
  Implementation, Nonconvex Optimization and its Applications, pp. 155--210.
  Springer (2006)

\bibitem{huang2020span}
Huang, X., Liang, X., Liu, Z., Li, L., Yu, Y., Li, Y.: Span: A stochastic
  projected approximate newton method.
\newblock Proceedings of the AAAI Conference on Artificial Intelligence
  \textbf{34}(02), 1520--1527 (2020)

\bibitem{hutchinson1989stochastic}
Hutchinson, M.F.: A stochastic estimator of the trace of the influence matrix
  for {L}aplacian smoothing splines.
\newblock Communications in Statistics-Simulation and Computation
  \textbf{18}(3), 1059--1076 (1989)

\bibitem{innes2019differentiable}
Innes, M., Edelman, A., Fischer, K., Rackauckas, C., Saba, E., Shah, V.B.,
  Tebbutt, W.: A differentiable programming system to bridge machine learning
  and scientific computing.
\newblock arXiv preprint arXiv:1907.07587  (2019)

\bibitem{Flux.jl-2018}
Innes, M., Saba, E., Fischer, K., Gandhi, D., Rudilosso, M.C., Joy, N.M.,
  Karmali, T., Pal, A., Shah, V.: Fashionable modelling with {F}lux.
\newblock CoRR \textbf{abs/1811.01457} (2018).
\newblock \urlprefix\url{https://arxiv.org/abs/1811.01457}

\bibitem{kelly2011new}
Kelly~Jr, J.L.: A new interpretation of information rate.
\newblock In: The Kelly capital growth investment criterion: theory and
  practice, pp. 25--34. World Scientific (2011)

\bibitem{kiwiel1990proximity}
Kiwiel, K.C.: Proximity control in bundle methods for convex nondifferentiable
  minimization.
\newblock Mathematical programming \textbf{46}(1), 105--122 (1990)

\bibitem{kiwiel2000efficiency}
Kiwiel, K.C.: Efficiency of proximal bundle methods.
\newblock Journal of Optimization Theory and Applications \textbf{104}(3),
  589--603 (2000)

\bibitem{lee2019inexact}
Lee, C.P., Wright, S.J.: Inexact successive quadratic approximation for
  regularized optimization.
\newblock Computational Optimization and Applications \textbf{72}(3), 641--674
  (2019)

\bibitem{lee2014proximal}
Lee, J.D., Sun, Y., Saunders, M.A.: Proximal {N}ewton-type methods for
  minimizing composite functions.
\newblock SIAM Journal on Optimization \textbf{24}(3), 1420--1443 (2014)

\bibitem{lemarechal1978nonsmooth}
Lemar{\'e}chal, C.: Nonsmooth optimization and descent methods.
\newblock RR-78-004 (1978)

\bibitem{lemarechal1995new}
Lemar{\'e}chal, C., Nemirovskii, A., Nesterov, Y.: New variants of bundle
  methods.
\newblock Mathematical programming \textbf{69}(1), 111--147 (1995)

\bibitem{lemarechal1997variable}
Lemar{\'e}chal, C., Sagastiz{\'a}bal, C.: Variable metric bundle methods: from
  conceptual to implementable forms.
\newblock Mathematical Programming \textbf{76}(3), 393--410 (1997)

\bibitem{levitin1966constrained}
Levitin, E.S., Polyak, B.T.: Constrained minimization methods.
\newblock USSR Computational mathematics and mathematical physics
  \textbf{6}(5), 1--50 (1966)

\bibitem{li2017inexact}
Li, J., Andersen, M.S., Vandenberghe, L.: Inexact proximal {N}ewton methods for
  self-concordant functions.
\newblock Mathematical Methods of Operations Research \textbf{85}(1), 19--41
  (2017)

\bibitem{lukvsan1998bundle}
Luk{\v{s}}an, L., Vl{\v{c}}ek, J.: A bundle-{N}ewton method for nonsmooth
  unconstrained minimization.
\newblock Mathematical Programming \textbf{83}(1), 373--391 (1998)

\bibitem{maclaurin2015autograd}
Maclaurin, D., Duvenaud, D., Adams, R.P.: Autograd: Effortless gradients in
  numpy.
\newblock In: ICML 2015 AutoML Workshop, vol. 238, p.~5 (2015)

\bibitem{marsh2007goodness}
Marsh, P.: Goodness of fit tests via exponential series density estimation.
\newblock Computational statistics \& data analysis \textbf{51}(5), 2428--2441
  (2007)

\bibitem{meyer2021hutch++}
Meyer, R.A., Musco, C., Musco, C., Woodruff, D.P.: Hutch++: Optimal stochastic
  trace estimation.
\newblock In: Symposium on Simplicity in Algorithms (SOSA), pp. 142--155. SIAM
  (2021)

\bibitem{mifflin1996quasi}
Mifflin, R.: A quasi-second-order proximal bundle algorithm.
\newblock Mathematical Programming \textbf{73}(1), 51--72 (1996)

\bibitem{mifflin1998quasi}
Mifflin, R., Sun, D., Qi, L.: Quasi-{N}ewton bundle-type methods for
  nondifferentiable convex optimization.
\newblock SIAM Journal on Optimization \textbf{8}(2), 583--603 (1998)

\bibitem{mordukhovich2020globally}
Mordukhovich, B.S., Yuan, X., Zeng, S., Zhang, J.: A globally convergent
  proximal {N}ewton-type method in nonsmooth convex optimization.
\newblock arXiv preprint arXiv:2011.08166  (2020)

\bibitem{nesterov2013gradient}
Nesterov, Y.: Gradient methods for minimizing composite functions.
\newblock Mathematical Programming \textbf{140}(1), 125--161 (2013)

\bibitem{nolan1953analytical}
Nolan, J.F.: Analytical differentiation on a digital computer.
\newblock Ph.D. thesis, Massachusetts Institute of Technology (1953)

\bibitem{noll2013bundle}
Noll, D.: Bundle method for non-convex minimization with inexact subgradients
  and function values.
\newblock In: Computational and analytical mathematics, pp. 555--592. Springer
  (2013)

\bibitem{ocpb:16}
O'Donoghue, B., Chu, E., Parikh, N., Boyd, S.: Conic optimization via operator
  splitting and homogeneous self-dual embedding.
\newblock Journal of Optimization Theory and Applications \textbf{169}(3),
  1042--1068 (2016).
\newblock \urlprefix\url{http://stanford.edu/~boyd/papers/scs.html}

\bibitem{scs}
O'Donoghue, B., Chu, E., Parikh, N., Boyd, S.: {SCS}: Splitting conic solver,
  version 2.1.2.
\newblock \url{https://github.com/cvxgrp/scs} (2019)

\bibitem{de2014convex}
de~Oliveira, W., Sagastiz{\'a}bal, C., Lemar{\'e}chal, C.: Convex proximal
  bundle methods in depth: a unified analysis for inexact oracles.
\newblock Mathematical Programming \textbf{148}(1), 241--277 (2014)

\bibitem{NEURIPS2019_9015}
Paszke, A., Gross, S., Massa, F., Lerer, A., Bradbury, J., Chanan, G., Killeen,
  T., Lin, Z., Gimelshein, N., Antiga, L., Desmaison, A., Kopf, A., Yang, E.,
  DeVito, Z., Raison, M., Tejani, A., Chilamkurthy, S., Steiner, B., Fang, L.,
  Bai, J., Chintala, S.: Py{T}orch: An imperative style, high-performance deep
  learning library.
\newblock Advances in Neural Information Processing Systems 32 pp. 8024--8035
  (2019).
\newblock
  \urlprefix\url{http://papers.neurips.cc/paper/9015-pytorch-an-imperative-style-high-performance-deep-learning-library.pdf}

\bibitem{rockafellar2002conditional}
Rockafellar, R.T., Uryasev, S.: Conditional value-at-risk for general loss
  distributions.
\newblock Journal of banking \& finance \textbf{26}(7), 1443--1471 (2002)

\bibitem{scheinberg2016practical}
Scheinberg, K., Tang, X.: Practical inexact proximal quasi-{N}ewton method with
  global complexity analysis.
\newblock Mathematical Programming \textbf{160}(1), 495--529 (2016)

\bibitem{schmidt2009optimizing}
Schmidt, M., Berg, E., Friedlander, M., Murphy, K.: Optimizing costly functions
  with simple constraints: A limited-memory projected quasi-{N}ewton algorithm.
\newblock In: Artificial Intelligence and Statistics, pp. 456--463. PMLR (2009)

\bibitem{schramm1992version}
Schramm, H., Zowe, J.: A version of the bundle idea for minimizing a nonsmooth
  function: Conceptual idea, convergence analysis, numerical results.
\newblock SIAM journal on optimization \textbf{2}(1), 121--152 (1992)

\bibitem{sra2012optimization}
Sra, S., Nowozin, S., Wright, S.J.: Optimization for machine learning.
\newblock MIT Press (2012)

\bibitem{teo2010bundle}
Teo, C.H., Vishwanathan, S.V.N., Smola, A., Le, Q.V.: Bundle methods for
  regularized risk minimization.
\newblock Journal of Machine Learning Research \textbf{11}(1) (2010)

\bibitem{cvxjl}
Udell, M., Mohan, K., Zeng, D., Hong, J., Diamond, S., Boyd, S.: Convex
  optimization in {J}ulia.
\newblock In: Proceedings of the Workshop for High Performance Technical
  Computing in Dynamic Languages, pp. 18--28 (2014)

\bibitem{CVAR_uryasev_rockafellar}
Uryasev, S., Rockafellar, R.T.: Conditional Value-at-Risk: Optimization
  Approach, pp. 411--435.
\newblock Springer US, Boston, MA (2001)

\bibitem{van2018incremental}
Van~Ackooij, W., Frangioni, A.: Incremental bundle methods using upper models.
\newblock SIAM Journal on Optimization \textbf{28}(1), 379--410 (2018)

\bibitem{wainwright2008graphical}
Wainwright, M., Jordan, M.: Graphical models, exponential families, and
  variational inference.
\newblock Now Publishers Inc (2008)

\bibitem{wu2010exponential}
Wu, X.: Exponential series estimator of multivariate densities.
\newblock Journal of Econometrics \textbf{156}(2), 354--366 (2010)

\bibitem{pmlr-v70-ye17a}
Ye, H., Luo, L., Zhang, Z.: Approximate {N}ewton methods and their local
  convergence.
\newblock In: D.~Precup, Y.W. Teh (eds.) Proceedings of the 34th International
  Conference on Machine Learning, \emph{Proceedings of Machine Learning
  Research}, vol.~70, pp. 3931--3939 (2017)

\bibitem{yu2010quasi}
Yu, J., Vishwanathan, S.V.N., G{\"u}nter, S., Schraudolph, N.N.: A
  quasi-{N}ewton approach to nonsmooth convex optimization problems in machine
  learning.
\newblock The Journal of Machine Learning Research \textbf{11}, 1145--1200
  (2010)

\bibitem{yue2019family}
Yue, M.C., Zhou, Z., So, A.M.C.: A family of inexact {SQA} methods for
  non-smooth convex minimization with provable convergence guarantees based on
  the {L}uo--{T}seng error bound property.
\newblock Mathematical Programming \textbf{174}(1), 327--358 (2019)

\end{thebibliography}

\end{document}